\documentclass[12pt]{article}
\usepackage{amssymb, amsmath, amsthm, amscd}
\usepackage[dvips]{graphics}
\usepackage[utf8]{inputenc}
\usepackage[all,cmtip]{xy}
\usepackage{bbm}
\usepackage{enumitem}
\usepackage{setspace}
\usepackage[colorlinks=true]{hyperref}
\hypersetup{colorlinks   = true}
\hypersetup{linkcolor=blue}
\begin{document}

\newtheorem{The}{Theorem}[section]
\newtheorem{Lem}[The]{Lemma}
\newtheorem{Prop}[The]{Proposition}
\newtheorem{Cor}[The]{Corollary}
\newtheorem{Rem}[The]{Remark}
\newtheorem{Obs}[The]{Observation}
\newtheorem{SConj}[The]{Standard Conjecture}
\newtheorem{Titre}[The]{\!\!\!\! }
\newtheorem{Conj}[The]{Conjecture}
\newtheorem{Question}[The]{Question}
\newtheorem{Prob}[The]{Problem}
\newtheorem{Def}[The]{Definition}
\newtheorem{Not}[The]{Notation}
\newtheorem{Claim}[The]{Claim}
\newtheorem{Conc}[The]{Conclusion}
\newtheorem{Ex}[The]{Example}
\newtheorem{Fact}[The]{Fact}
\newtheorem{Formula}[The]{Formula}
\newtheorem{Formulae}[The]{Formulae}
\newcommand{\C}{\mathbb{C}}
\newcommand{\R}{\mathbb{R}}
\newcommand{\N}{\mathbb{N}}
\newcommand{\Z}{\mathbb{Z}}
\newcommand{\Q}{\mathbb{Q}}
\newcommand{\Proj}{\mathbb{P}}
\newcommand{\Rc}{\mathcal{R}}
\newcommand{\Oc}{\mathcal{O}}
\newcommand\houda[1]{{\textcolor{green}{#1}}}
\newcommand\dan[1]{{\textcolor{blue}{#1}}}
\begin{center}

{\Large\bf Positivity Cones under Deformations of Complex Structures}

\end{center}

\begin{center}

{\large Houda Bellitir and Dan Popovici}

\end{center}

\vspace{1ex}

\noindent{\small{\bf Abstract.} We investigate connections between the sGG property of compact complex manifolds, defined in earlier work by the second author and L. Ugarte by the requirement that every Gauduchon metric be strongly Gauduchon, and a possible degeneration of the Fr\"olicher spectral sequence. In the first approach that we propose, we prove a partial degeneration at $E_2$ and we introduce a positivity cone in the $E_2$-cohomology of bidegree $(n-2,\,n)$ of the manifold that we then prove to behave lower semicontinuously under deformations of the complex structure. In the second approach that we propose, we introduce an analogue of the $\partial\bar\partial$-lemma property of compact complex manifolds for any real non-zero constant $h$ using the partial twisting $d_h$, introduced recently by the second author, of the standard Poincar\'e differential $d$. We then show, among other things, that this $h$-$\partial\bar\partial$-property is deformation open.}

\vspace{1ex}

\section{Introduction}\label{section:introduction}

In this paper, we explore a few connections between the metric geometry of compact complex manifolds and some Hodge-theoretic aspects related to the Fr\"olicher spectral sequence (a classical object linking the differential and the complex structures, see e.g. a reminder of the definition in $\S.$\ref{section:preliminaries} below) of these manifolds. 

It is well known that the existence of a K\"ahler metric implies the best possible degeneration (i.e. at the first page $E_1$) of this spectral sequence. However, compact complex manifolds admit only rarely K\"ahler metrics and no other metric property is currently known to imply the degeneration at some page of this spectral sequence. The following conjecture was proposed and solved in a special case in [Pop16].

\begin{Conj}\label{Conj:SKT_E2} Let $X$ be a compact complex manifold. If an SKT metric (i.e. a $C^\infty$ positive definite $(1,\,1)$-form $\omega$ such that $\partial\bar\partial\omega=0$) exists on $X$, the Fr\"olicher spectral sequence of $X$ degenerates at the second page $E_2$.

\end{Conj}

The general case of this conjecture, that will certainly have a role to play in the classification theory of compact complex manifolds, remains open. In this paper, we investigate a possible variant of it in which SKT metrics are replaced by strongly Gauduchon (sG) metrics and the degeneration issue is often confined to the cohomology of $X$ in a degree close to the maximal one. We take our cue from a result of Ceballos-Otal-Ugarte-Villacampa [COUV16, Theorem 5.6] asserting that {\it the Fr\"olicher spectral sequence of any $6$-real-dimensional nilmanifold endowed with an invariant complex structure and carrying an sG metric degenerates at $E_2$} and from a possible generalisation of this statement they wondered about in 

\begin{Question}([COUV16, Question 5.7])\label{Question:COUV_sG-F} Does the Fr\"olicher spectral sequence of any $3$-dimensional compact complex manifold carrying an sG metric degenerate at $E_2$? \end{Question}

\vspace{2ex} 

Let $X$ be a compact complex manifold with $\mbox{dim}_\C X=n$. Recall that Hermitian metrics, defined as $C^\infty$ positive definite $(1,\,1)$-forms $\omega$, always exist on $X$. Even {\it Gauduchon metrics} $\omega$, defined as Hermitian metrics satisfying the extra condition $\partial\bar\partial\omega^{n-1}=0$, always exist (cf. [Gau77]). However, {\it strongly Gauduchon (sG) metrics}, introduced in [Pop13] in the context of deformations of complex structures and defined by the stronger requirement that $\partial\omega^{n-1}$ be $\bar\partial$-exact, need not exist. Compact complex manifolds that admit sG metrics were called {\it strongly Gauduchon (sG) manifolds} and the notion covers a wide range of manifolds and considerably enlarges the class of compact K\"ahler manifolds and their bimeromorphic models ($=$ the so-called {\it Fujiki class ${\cal C}$ manifolds}).   

Section \ref{section:E_2sG_cone} of this paper investigates elements of the $E_2$ degeneration of the Fr\"olicher spectral sequence of {\it sGG manifolds}, a class of compact complex manifolds introduced and studied in [PU18a]. They are defined by the requirement that every Gauduchon metric be strongly Gauduchon. In particular, sG metrics exist on every sGG manifold. Moreover, thanks to [PU18a, Lemma 1.3], an $n$-dimensional compact complex manifold $X$ is sGG if and only if the following special case of the $\partial\bar\partial$-lemma holds on $X$: \\

{\it for every $d$-closed $(n,\,n-1)$-form $\Gamma$ on $X$, if $\Gamma$ is $\partial$-exact, then $\Gamma$ is also $\bar\partial$-exact.}  \\

We exploit this fact in at least two ways in section \ref{section:E_2sG_cone}: \\

$(1)$\, by showing that every sGG manifold has the partial Fr\"olicher $E_2$ degeneration property that the map $d_2$ acting in bidegree $(n-2,\,n)$ on the second page of its Fr\"olicher spectral sequence vanishes identically (cf. Proposition \ref{Prop:sGG_E2});

\vspace{1ex}

$(2)$\, by introducing three versions ${\cal S}_X\subset E_2^{n-2,\,n}(X)$, $\widetilde{\cal S}_X\subset H_{DR}^{2n-2}(X,\,\R)$ (cf. Definition \ref{Def:cone_S}) and $\widehat{\cal S}_X\subset E_2^{n-2,\,n}(X)$ (cf. (\ref{eqn:E2sG-cone-hat})) of a positivity cone that we call the {\it $E_2sG$-cone} of a given compact complex $n$-dimensional manifold $X$. The term $E_2sG$ refers to the fact that this cone consists of (double) cohomology classes of bidegree $(n-2,\,n)$ that arise on the second page of the Fr\"olicher spectral sequence of $X$ and are thus a refinement of the Dolbeault cohomology.  

The surprising aspect is that we thus effectively introduce a notion of positivity for cohomology classes of bidegree $(n-2,\,n)$ that runs counter to the familiar notions of positivity that exist in all the bidegrees $(p,\,p)$ (but not $(p,\,q)$ with $p\neq q$) in complex geometry. This is done by using strongly Gauduchon metrics that enable the existing notion of positivity in bidegree $(n-1,\,n-1)$ to carry over to the bidegree $(n-2,\,n)$ in a natural way. The $E_2sG$-cone is empty if the manifold $X$ is not strongly Gauduchon. It depends on the complex structure of $X$. We study this dependence by proving the following

\begin{The}\label{The:E_2sG_semicontinuity_introd} Let $\pi:{\cal X}\longrightarrow\Delta$ be a holomorphic family of compact complex $n$-dimensional manifolds over a ball $\Delta\subset\C^N$ centred at the origin. Suppose that the fibre $X_0:=\pi^{-1}(0)$ is an {\bf sGG manifold}. 

Then, the De Rham $E_2sG$-cone $\widetilde{\cal S}_{X_t}\subset H_{DR}^{2n-2}(X,\,\R)$ of the fibre $X_t:=\pi^{-1}(t)\subset{\cal X}$ varies in a {\bf lower semicontinuous} way with $t\in\Delta$ varying in a small enough neighbourhood of $0\in\Delta$. 

\end{The}

\vspace{2ex}

By $X$ we mean the smooth manifold that underlies all the fibres $X_t$ of the family (known to be $C^\infty$ trivial by the classical Ehresmann Theorem, but in general not holomorphically trivial), so the real De Rham cohomology group $H_{DR}^{2n-2}(X,\,\R)$ of degree $2n-2$ is independent of the (complex structure of the) fibre $X_t$. The meaning of {\it lower semicontinuous} in connection with the dependence on $t$ of the De Rham $E_2sG$-cone $\widetilde{\cal S}_{X_t}$ is made precise in Theorem \ref{The:E_2sG_semicontinuity}.

We also use the Serre-type duality (proved in the forthcoming joint paper [PU18b] of the second author with L. Ugarte by means of the pseudo-differential Laplacian introduced in [Pop16]) between any pair $(E_2^{p,\,q}(X),\, E_2^{n-p,\,n-q}(X))$ of vector spaces of complementary bidegrees featuring on the second page of the Fr\"olicher spectral sequence of $X$. In this way, we prove the following analogue in our context of Lamari's duality [Lam99, Lemma 3.3] between the pseudo-effective cone and the closure of the Gauduchon cone of any compact complex manifold.

\begin{Prop}\label{Prop:cone-duality_introd} Let $X$ be an $n$-dimensional {\bf sGG} compact complex manifold on which an arbitrary Hermitian metric $\gamma$ has been fixed. The dual of the closure of the cone \begin{equation}\label{eqn:E2sG-cone-hat_introd}\widehat{\cal S}_X = \bigg\{\bigg[[\Gamma^{n-2,\,n}]_{\bar\partial}\bigg]_{d_1}\,\mid\, \exists\,\omega \hspace{1ex} \mbox{Hermitian metric such that} \hspace{1ex} \partial\Gamma^{n-2,\,n} = -\bar\partial\omega^{n-1} \bigg\}\subset E_2^{n-2,\,n}(X)\end{equation}

\noindent under the duality $E_2^{2,\,0}(X)\times E_2^{n-2,\,n}(X)\longrightarrow\C$
 is the closed convex cone in $E_2^{2,\,0}(X)$ consisting of the $E_2$-classes $[[\theta^{2,\,0}]_{\bar\partial}]_{d_1}$ ``representable'' by $\gamma$-positive currents $\tau^{2,\,0} : C^\infty_{n-2,\,n}(X,\,\R)_{\gamma}\longrightarrow\R$.

\end{Prop}

We refer to Proposition \ref{Prop:cone-duality} for a more precise statement and to Definition \ref{Def:positivity_skewed-bideg} for the notion of {\it $\gamma$-positive current of bidegree $(2,\,0)$}. This is another extension of the classical notions of positivity in complex geometry to a bidegree $(p,\,q)$ with $p\neq q$. The bidegree $(2,\,0)$ is especially important on {\it holomorphic symplectic (not necessarily K\"ahler) manifolds} and this lead will hopefully be investigated in future work.  

\vspace{3ex}

Section \ref{section:h-dd-bar} of this paper takes up the Fr\"olicher degeneration issue and the variation of the complex structure from a different point of view that still ties in with the theory of sGG manifolds. 

 Let $X$ be a compact complex manifold with $\mbox{dim}_\C X=n$. In [Pop17], for every positive constant $h$, the differential operator $d=\partial + \bar\partial$ associated with the smooth structure of $X$ was modified to \begin{equation}\label{eqn:d_h_def}d_h:=h\partial + \bar\partial:C^\infty_k(X,\,\C)\to C^\infty_{k+1}(X,\,\C), \hspace{6ex} k\in\{0,\dots , 2n\},\footnote{Throughout the paper, $C^\infty_k(X,\,\C)$ will stand for the space of $\C$-valued $C^\infty$ $k$-forms on $X$.}\end{equation} \noindent by rescaling its $(1,\,0)$-part in the splitting induced by the complex structure of $X$. Unlike $d$, the operators $d_h$ depend on the complex structure of $X$ while also sharing some properties with $d$. The most striking of these is that the $d_h$-cohomology of $X$, relying on the integrability property $d_h^2=0$ and defined for every $k\in\{0,\dots , 2n\}$ by $$H^k_{d_h}(X,\,\C)=\ker(d_h:C^\infty_k(X,\,\C)\to C^\infty_{k+1}(X,\,\C))\bigg\slash\mbox{Im}\,(d_h:C^\infty_{k-1}(X,\,\C)\to C^\infty_k(X,\,\C)),$$ \noindent is isomorphic to the De Rham cohomology of $X$ via the isomorphism $H^k_{DR}(X,\,\C)\ni\{u\}_{DR}\mapsto\{\theta_h u\}_{d_h}\in H^k_{d_h}(X,\,\C)$ induced by the pointwise isomorphism

$$\theta_h:\Lambda^{p,\,q}T^\star X\to \Lambda^{p,\,q}T^\star X, \hspace{3ex} u\mapsto\theta_h u:=h^p\,u,$$

\noindent at the level of pure-type differential forms on $X$. The operators $d_h$ capture in a certain sense the relationships between the smooth and the complex structures of $X$.

On the other hand, there exist in the literarure at least two notions of the $\partial\bar\partial$-lemma being satisfied by a compact complex manifold that capture the smooth structure-complex structure relationship. 

\vspace{1ex}

(a)\, An early notion that has been used by many authors appeared in [DGMS75, Lemmas 5.11 and 5.15]. It defines the fulfilment of the $\partial\bar\partial$-lemma (or of the equivalent $dd^c$-lemma) on a compact complex manifold by the requirement that any smooth differential form $u$ of any degree (but {\it not necessarily of pure type}) that is both $\partial$-closed and $\bar\partial$-closed satisfies the following implication:

\begin{equation}\label{eqn:d-exact_implies_dd-bar-exact_introd}u\in\mbox{Im}\,d \implies u\in\mbox{Im}\,(\partial\bar\partial).\end{equation}

\noindent This implication is actually an equivalence since the reverse implication holds trivially on any manifold. 

For example, this is what was meant by the $\partial\bar\partial$-lemma holding on a given manifold in [AT12]. 

\vspace{1ex}

(b)\, In [Pop14], the term {\it $\partial\bar\partial$-manifold} was introduced to mean that a given compact complex manifold $X$ satisfies the $\partial\bar\partial$-lemma if for any $d$-closed {\it pure-type} form $u$ on $X$ the following exactness properties are equivalent:

\begin{equation}\label{eqn:all-exactness-equiv_introd}u\in\mbox{Im}\,d \iff u\in\mbox{Im}\,\partial  \iff u\in\mbox{Im}\,\bar\partial \iff u\in\mbox{Im}\,(\partial\bar\partial).\end{equation} 

\noindent The last property trivially implies the others, so the above equivalences reduce to each of the other three forms of exactness implying $(\partial\bar\partial)$-exactness. Since $u$ is of pure type, the $d$-closedness assumption on $u$ is equivalent to $u$ being assumed both $\partial$-closed and $\bar\partial$-closed. 

On the face of it, condition (\ref{eqn:all-exactness-equiv_introd}) is more restrictive than (\ref{eqn:d-exact_implies_dd-bar-exact_introd}), but (\ref{eqn:all-exactness-equiv_introd}) is only required to apply to {\it pure-type forms}. It is implicit in [DGMS75, Lemma 5.15] that version (a) of the $\partial\bar\partial$-lemma condition (required to hold on all, not necessarily pure-type forms) actually implies version (b).

For every constant $h\in\R\setminus\{0\}$, we introduce in this paper the notion of {\it $h$-$\partial\bar\partial$-manifold} that implies the above version (a) (hence also version (b)) of the $\partial\bar\partial$-condition. The idea is to use the operators $d_h$ to capture the interplay between the smooth structure and the complex structure.

\begin{Def}\label{Def:h-dd-bar-manifolds_introd} Let $h\in\R\setminus\{0\}$ be an arbitrary constant. A compact complex manifold $X$ with $\mbox{dim}_\C X=n$ is said to be an {\bf h-$\partial\bar\partial$-manifold} if for every $k\in\{0,1,\dots ,2n\}$ and every $k$-form $u\in\ker d_h\cap\ker d_{-h^{-1}}$, the following exactness conditions are equivalent:

$$u\in\mbox{Im}\,d_h \iff u\in\mbox{Im}\,d_{-h^{-1}} \iff u\in\mbox{Im}\,d \iff u\in\mbox{Im}\,(d_h\,d_{-h^{-1}}) = \mbox{Im}\,(\partial\bar\partial).$$

\end{Def}

We prove in Corollary \ref{Cor:equivalence_h-dd-bar_properties} that an equivalent property is obtained by the removal of the condition $u\in\mbox{Im}\,d$ from the above sequence of equivalences. Note that the forms $u$ are not required to be of pure type in Definition \ref{Def:h-dd-bar-manifolds_introd}. Unlike $\partial$, $\bar\partial$ and $\partial\bar\partial$, the operators $d_h$ do not map pure-type forms to pure-type forms, so the proof of the implication ``$(a)_n\Longrightarrow(b)_n$ `` in [DGMS75, Lemma 5.15] does not seem to adapt easily to yield a proof of the possible implication ``version (a) of the $\partial\bar\partial$-property $\Longrightarrow$ the $h$-$\partial\bar\partial$-property``. Actually, we do not know whether this last implication holds. A priori, the $h$-$\partial\bar\partial$-property is stronger when $h\notin\{-1,\,1\}$. 

Introducing $d_h$-analogues of the standard Bott-Chern and Aeppli cohomologies and following the pattern of Wu's proof in [Wu06] of the stability under small deformations of the complex structure of the standard $\partial\bar\partial$-property, we prove that the analogous statement holds for our $h$-$\partial\bar\partial$-property.

\begin{The}\label{The:h-dd-bar-openness_introd} Fix an arbitrary constant $h\in\R\setminus\{0\}$. The $h$-$\partial\bar\partial$-property of compact complex manifolds is open under deformations of the complex structure.

\end{The}

See Theorem \ref{The:h-dd-bar-openness} for a more precise statement.

One last explanation is in order about the choice of pairing $d_h$ with $d_{-h^{-1}}$, rather than with the more natural-looking $\overline{d}_h= h\,d_{h^{-1}}$. This choice is forced on us by a formula of the Bochner-Kodaira-Nakano type (the {\bf h-BKN identity}) that we establish in Theorem \ref{Thm:BKN} as a consequence of what we call {\it $h$-commutation relations} that we compute for the operators $d_h$ and for an arbitrary Hermitian metric in Lemma \ref{Lem:h-commutation-rel}. The pair $(d_h,\,d_{-h^{-1}})$ generalises the classical pair $(d,\,d^c)$ since for $h=-1$, $d_{-1}$ is a constant multiple (which does not change either the kernel or the image) of $d^c = -JdJ = i(\bar\partial-\partial) = i\,d_{-1}$.

\vspace{3ex}

\noindent {\bf Acknowledgments.} This paper is part of the first author's PhD thesis under the co-supervision of Sa\"{\i}d Asserda and the second author. The first author is grateful to Professor Sa\"{\i}d Asserda for his constant support, patient guiding of her work throughout the duration of her thesis and careful reading of the paper. Both authors are grateful to Ahmed Zeriahi for fruitful discussions and for his interest in this work.

\section{Preliminaries}\label{section:preliminaries}

Let $X$ be a compact complex manifold with $\mbox{dim}_\C X=n$. 

Let $(E_r,\,d_r)$ be the Fr\"olicher spectral sequence of $X$. It relates the De Rham cohomology of $X$ to its Dolbeault cohomology in the following way. The $0^{th}$ page features the (infinite-dimensional) $\C$-vector spaces $E_0^{p,\,q}=C^\infty_{p,\,q}(X,\,\C)$ of $C^\infty$ forms of arbitrary bidegree $(p,\,q)$ with $0\leq p,q\leq n$ and the linear maps

$$\dots \stackrel{d_0}{\longrightarrow}E_0^{p,\,q}(X)\stackrel{d_0}{\longrightarrow} E_0^{p,\,q+1}(X)\stackrel{d_0}{\longrightarrow}\dots,  \hspace{3ex} \mbox{where} \hspace{1ex} d_0 =\bar\partial.$$

 The $1^{st}$ page is defined as the cohomology of the $0^{th}$ page and consists of the (finite-dimensional) Dolbeault cohomology groups $$E_1^{p,\,q}(X)=H^{p,\,q}(X,\,\C) = \ker(d_0:E_0^{p,\,q}(X)\to E_0^{p,\,q+1}(X))/\mbox{Im}\,(d_0:E_0^{p,\,q-1}(X)\to E_0^{p,\,q}(X))$$

\noindent and the linear maps $$\dots \stackrel{d_1}{\longrightarrow}E_1^{p,\,q}(X)\stackrel{d_1}{\longrightarrow} E_1^{p+1,\,q}(X)\stackrel{d_1}{\longrightarrow}\dots$$ \noindent defined by $\partial$ in cohomology in the following way: $d_1([\alpha]_{\bar\partial})=[\partial\alpha]_{\bar\partial}$ for all $[\alpha]_{\bar\partial}\in E_1^{p,\,q}$. 

We then continue by induction and define the $r^{th}$ page as the cohomology of the $(r-1)^{st}$ page, namely $$E_r^{p,\,q}(X) = \ker(d_{r-1}:E_{r-1}^{p,\,q}(X)\to E_{r-1}^{p+r-1,\,q-r+2}(X))/\mbox{Im}\,(d_{r-1}:E_{r-1}^{p-r+1,\,q+r-2}(X)\to E_{r-1}^{p,\,q}(X)),$$ \noindent where the linear maps $d_r$ on each page $r$ are of bidegree $(r,\,-r+1)$. In particular, the $2^{nd}$ page is of the form $$\dots \stackrel{d_2}{\longrightarrow}E_2^{p,\,q}(X)\stackrel{d_2}{\longrightarrow} E_2^{p+2,\,q-1}(X)\stackrel{d_2}{\longrightarrow}\dots,$$ \noindent where each space $E_2^{p,\,q}(X)$ consists of double cohomology classes $\bigg[[\alpha]_{\bar\partial}\bigg]_{d_1}$ of smooth $(p,\,q)$-forms $\alpha$ satisfying the conditions

\begin{equation}\label{eqn:E_2_rep}\bar\partial\alpha=0 \hspace{3ex} \mbox{and} \hspace{3ex} \partial\alpha\in\mbox{Im}\,\bar\partial,\end{equation}

\noindent while each map $d_2:E_2^{p,\,q}(X)\longrightarrow E_2^{p+2,\,q-1}(X)$ is defined by

\begin{equation}\label{eqn:E_2_map_def}d_2\bigg(\bigg[[\alpha]_{\bar\partial}\bigg]_{d_1}\bigg) = \bigg[[\partial u_1]_{\bar\partial}\bigg]_{d_1} \hspace{3ex} \mbox{where} \hspace{1ex} u_1 \hspace{1ex} \mbox{is any form such that} \hspace{1ex} \partial\alpha = \bar\partial u_1.\end{equation}

\noindent The definition of $d_2$ is independent of the choice of the $\bar\partial$-potential $u_1$. 

It is easy to check that the vanishing condition for an arbitrary element $\bigg[[\alpha]_{\bar\partial}\bigg]_{d_1}\in E_2^{p,\,q}(X)$ is

\begin{equation}\label{eqn:E_2_class-vanishing}\bigg[[\alpha]_{\bar\partial}\bigg]_{d_1} = 0 \iff \exists u\in C^\infty_{p-1,\,q}(X,\,\C)\cap\ker\bar\partial \hspace{1ex} \mbox{and} \hspace{1ex}  v\in C^\infty_{p,\,q-1}(X,\,\C) \hspace{1ex} \mbox{such that} \hspace{1ex} \alpha = \partial u + \bar\partial v.\end{equation}

\section{The $E_2sG$ cone}\label{section:E_2sG_cone}

Let $X$ be a compact complex manifold with $\mbox{dim}_\C X=n$.

\subsection{Complex-valued cohomology}\label{subsection:C-valued}

 We shall first deal essentially with the De Rham cohomology of $X$ with values in $\C$ and the standard space $E_2^{n-2,\,n}(X)$ on the second page of the Fr\"olicher spectral sequence.

\begin{Prop}\label{Prop:T_def-image} The following canonical linear map

$$T: H^{2n-2}_{DR}(X,\,\C)\longrightarrow E_2^{n-2,\,n}(X), \hspace{3ex}  \{\alpha\}_{DR}\mapsto\bigg[[\alpha^{n-2,\,n}]_{\bar\partial}\bigg]_{d_1},$$

\noindent is well defined and its image is given by

\begin{equation}\label{eqn:Im_T}\mbox{Im}\,T =\ker d_2^{n-2,\,n},  \end{equation}

\noindent where $d_2^{n-2,\,n}: E_2^{n-2,\,n}(X)\to E_2^{n,\,n-1}(X)$ is the $d_2$-map acting in bidegree $(n-2,\,n)$.

\end{Prop}

\noindent {\it Proof.} Let $\{\alpha\}_{DR}\in H^{2n-2}_{DR}(X,\,\C)$ be an arbitrary class and let $\alpha = \alpha^{n,\,n-2} + \alpha^{n-1,\,n-1} + \alpha^{n-2,\,n}$ be an arbitrary representative, where the $\alpha^{p,\,q}$'s are the components of $\alpha$ of types $(p,\,q)$. The condition $d\alpha =0$ is equivalent to

$$\partial\alpha^{n-1,\,n-1} + \bar\partial\alpha^{n,\,n-2} =0 \hspace{3ex} \mbox{and} \hspace{3ex} \bar\partial\alpha^{n-1,\,n-1} + \partial\alpha^{n-2,\,n} =0.$$

\noindent Since $\bar\partial\alpha^{n-2,\,n}=0$ and $\partial\alpha^{n-2,\,n} = -\bar\partial\alpha^{n-1,\,n-1}\in\mbox{Im}\,\bar\partial$, $\alpha^{n-2,\,n}$ defines a class $[[\alpha^{n-2,\,n}]_{\bar\partial}]_{d_1}$ in $E_2^{n-2,\,n}(X)$. 

To show well-definedness for $T$, we still have to show that the definition is independent of the choice of representative $\alpha$ of the De Rham class $\{\alpha\}_{DR}$. This is equivalent to showing that $T$ maps $0\in H^{2n-2}_{DR}(X,\,\C)$ to $0\in E_2^{n-2,\,n}(X)$. Let $\alpha\in C^\infty_{2n-2}(X,\,\C)$ be $d$-exact. Then, there exists $\beta=\beta^{n,\,n-3} + \beta^{n-1,\,n-2} + \beta^{n-2,\,n-1} + \beta^{n-3,\,n} $ a $(2n-3)$-form such that $\alpha=d\beta$. This amounts to

\begin{equation}\label{eqn:exactness-condition}\alpha^{n,\,n-2} = \partial\beta^{n-1,\,n-2} + \bar\partial\beta^{n,\,n-3}, \hspace{2ex} \alpha^{n-1,\,n-1} = \partial\beta^{n-2,\,n-1} + \bar\partial\beta^{n-1,\,n-2}  \hspace{2ex} \mbox{and} \hspace{2ex} \alpha^{n-2,\,n} = \partial\beta^{n-3,\,n} + \bar\partial\beta^{n-2,\,n-1}. \end{equation} 

\noindent Since $\bar\partial\beta^{n-3,\,n}=0$ for bidegree reasons, the last identity in (\ref{eqn:exactness-condition}) shows, thanks to (\ref{eqn:E_2_class-vanishing}), that $[[\alpha^{n-2,\,n}]_{\bar\partial}]_{d_1}=0$ in $E_2^{n-2,\,n}(X)$. 

Let us now prove the inclusion $\ker d_2^{n-2,\,n}\subset\mbox{Im}\,T$ in (\ref{eqn:Im_T}). Let $[[\alpha^{n-2,\,n}]_{\bar\partial}]_{d_1}\in E_2^{n-2,\,n}(X)$ such that $d_2([[\alpha^{n-2,\,n}]_{\bar\partial}]_{d_1})=0$. Thanks to (\ref{eqn:E_2_rep}), (\ref{eqn:E_2_map_def}) and (\ref{eqn:E_2_class-vanishing}), there exist forms $\Omega^{n-1,\,n-1}\in C^\infty_{n-1,\,n-1}(X,\,\C)$, $u\in C^\infty_{n-1,\,n-1}(X,\,\C)$ with $u\in\ker\bar\partial$ and $v\in C^\infty_{n,\,n-2}(X,\,\C)$ such that 

$$\partial\alpha^{n-2,\,n} = -\bar\partial\Omega^{n-1,\,n-1}  = -\bar\partial(\Omega^{n-1,\,n-1}-u)\hspace{2ex} \mbox{and} \hspace{2ex} \partial\Omega^{n-1,\,n-1} = \partial u + \bar\partial v.$$

\noindent If we put $\alpha:=\alpha^{n-2,\,n} + (\Omega^{n-1,\,n-1} -u) -v$, we see that 

$$d\alpha =0  \hspace{2ex} \mbox{and} \hspace{2ex} T(\{\alpha\}_{DR}) =[[\alpha^{n-2,\,n}]_{\bar\partial}]_{d_1}.$$

Let us now prove the reverse inclusion $\ker d_2^{n-2,\,n}\supset\mbox{Im}\,T$ in (\ref{eqn:Im_T}). Let $[[\alpha^{n-2,\,n}]_{\bar\partial}]_{d_1}\in\mbox{Im}\, T$. This means that $\alpha^{n-2,\,n}$ is the $(n-2,\,n)$-component of a $d$-closed $(2n-2)$-form $\alpha = \alpha^{n,\,n-2} + \alpha^{n-1,\,n-1} + \alpha^{n-2,\,n}$. As already noticed, the condition $d\alpha=0$ is equivalent to

$$\partial\alpha^{n-1,\,n-1} + \bar\partial\alpha^{n,\,n-2}=0 \hspace{2ex} \mbox{and} \hspace{2ex} \partial\alpha^{n-2,\,n} + \bar\partial\alpha^{n-1,\,n-1} =0.$$

\noindent On the other hand, $d_2([[\alpha^{n-2,\,n}]_{\bar\partial}]_{d_1}) = - [[\partial\alpha^{n-1,\,n-1}]_{\bar\partial}]_{d_1} =  [[\bar\partial\alpha^{n,\,n-2}]_{\bar\partial}]_{d_1} =0\in E_2^{n,\,n-1}$ since even $[\bar\partial\alpha^{n,\,n-2}]_{\bar\partial}=0\in H^{n,\,n-1}_{\bar\partial}(X,\,\C)$.  \hfill $\Box$

\vspace{2ex}

As a consequence, we obtain the following criterion for partial degeneration at $E_2$ of the Fr\"olicher spectral sequence of $X$.

\begin{Cor}\label{Cor:d_2_vanishing} The canonical map $T$ defined in Proposition \ref{Prop:T_def-image} is {\bf surjective} if and only if the map $d_2$ {\bf vanishes identically} in bidegree $(n-2,\,n)$.

\end{Cor}

We now show that the sGG assumption on the ambient manifold $X$ (see [PU18a] for the definition and a study of the class of sGG manifolds) suffices to guarantee the partial degeneration property mentioned above.

\begin{Prop}\label{Prop:sGG_E2} Let $X$ be a compact complex manifold with $\mbox{dim}_\C X=n$. If $X$ is {\bf sGG}, the map $d_2^{n-2,\,n}:E_2^{n-2,\,n}(X)\longrightarrow E_2^{n,\,n-1}(X)$ on the $2^{nd}$ page of the Fr\"olicher spectral sequence of $X$ {\bf vanishes identically} (equivalently, the canonical linear map $T:H^{2n-2}_{DR}(X,\,\C)\longrightarrow E_2^{n-2,\,n}(X)$ of Proposition \ref{Prop:T_def-image} is {\bf surjective}).

\end{Prop}

\noindent {\it Proof.} To prove that $T$ is surjective, let $[[\alpha^{n-2,\,n}]_{\bar\partial}]_{d_1}\in E_2^{n-2,\,n}(X)$ and let $\alpha^{n-2,\,n}$ be an arbitrary $(n-2,\,n)$-form representing this double class. Then $\partial\alpha^{n-2,\,n}$ is $\bar\partial$-exact, so there exists an $(n-1,\,n-1)$-form $\Omega^{n-1,\,n-1}$ such that 

\begin{equation}\label{eqn:choice-Omega}\partial\alpha^{n-2,\,n} = -\bar\partial\Omega^{n-1,\,n-1}.\end{equation} 

 Now, the sGG assumption on $X$ implies that $\bar\partial\overline{\Omega^{n-1,\,n-1}}$ is $\partial$-exact. Indeed, this is equivalent to $\partial\Omega^{n-1,\,n-1}$ being $\bar\partial$-exact. Meanwhile, $\partial\Omega^{n-1,\,n-1}$ is a $d$-closed and $\partial$-exact $(n,\,n-1)$-form, so (iii) of Lemma 1.3. in [PU18a] implies that $\partial\Omega^{n-1,\,n-1}$ is also $\bar\partial$-exact (thanks to $X$ being sGG).

Thus, there exists $\beta^{n-2,\,n}\in C^\infty_{n-2,\,n}(X,\,\C)$ such that

\begin{equation}\label{eqn:choice-beta}\partial\beta^{n-2,\,n} = -\bar\partial\overline{\Omega^{n-1,\,n-1}}.\end{equation}

\noindent Consequently, $\partial(\alpha^{n-2,\,n} + \beta^{n-2,\,n}) = -\bar\partial(\Omega^{n-1,\,n-1} + \overline{\Omega^{n-1,\,n-1}})$, hence 

$$\Gamma_1:= \overline{(\alpha^{n-2,\,n} + \beta^{n-2,\,n})} + (\Omega^{n-1,\,n-1} + \overline{\Omega^{n-1,\,n-1}}) + (\alpha^{n-2,\,n} + \beta^{n-2,\,n})$$

\noindent is a $(2n-2)$-form such that $d\Gamma_1=0$ and 

\begin{equation}\label{eqn:T_gamma_1}T(\{\Gamma_1\}_{DR}) = [[\alpha^{n-2,\,n} + \beta^{n-2,\,n}]_{\bar\partial}]_{d_1}\in E_2^{n-2,\,n}(X).\end{equation}

\noindent Note that $\beta^{n-2,\,n}$ defines indeed an $E_2$-class since it is $\bar\partial$-closed (for bidegree reasons) and $\partial\beta^{n-2,\,n}$ is $\bar\partial$-exact (by construction).

We also get $\partial(\alpha^{n-2,\,n} - \beta^{n-2,\,n}) = -\bar\partial(\Omega^{n-1,\,n-1} - \overline{\Omega^{n-1,\,n-1}})$, hence 

$$\Gamma_2:= (\overline{\beta^{n-2,\,n}} - \overline{\alpha^{n-2,\,n}}) + (\Omega^{n-1,\,n-1} - \overline{\Omega^{n-1,\,n-1}}) + (\alpha^{n-2,\,n} - \beta^{n-2,\,n})$$

\noindent is a $(2n-2)$-form such that $d\Gamma_2=0$ and 

\begin{equation}\label{eqn:T_gamma_2}T(\{\Gamma_2\}_{DR}) = [[\alpha^{n-2,\,n} - \beta^{n-2,\,n}]_{\bar\partial}]_{d_1}\in E_2^{n-2,\,n}(X).\end{equation}

Putting (\ref{eqn:T_gamma_1}) and (\ref{eqn:T_gamma_2}) together, we finally get
a $d$-closed $(2n-2)$-form 

$$\frac{\Gamma_1+\Gamma_2}{2} = \overline{\beta^{n-2,\,n}} + \Omega^{n-1,\,n-1} + \alpha^{n-2,\,n}$$ 

\noindent satisfying the condition

$$T\bigg(\bigg\{\frac{\Gamma_1+\Gamma_2}{2}\bigg\}_{DR}\bigg) = \bigg[[\alpha^{n-2,\,n}]_{\bar\partial}\bigg]_{d_1}\in E_2^{n-2,\,n}(X).$$

\noindent Thus, $T$ is surjective.  \hfill $\Box$

\vspace{3ex}

When the manifold $X$ is sGG, we can take the surjectivity of the canonical linear map $T:H^{2n-2}_{DR}(X,\,\C)\longrightarrow E_2^{n-2,\,n}(X)$ of Proposition \ref{Prop:T_def-image} further by showing that every Hermitian metric $\omega$ on $X$ defines a natural {\it injection} of $E_2^{n-2,\,n}(X)$ into $H^{2n-2}_{DR}(X,\,\C)$ that is a section of $T$. We will need the following Laplace-type pseudo-differential operator

$$\widetilde\Delta:=\partial p''\partial^\star + \partial^\star p''\partial + \bar\partial\bar\partial^\star + \bar\partial^\star\bar\partial : C^\infty_{p,\,q}(X,\,\C)\longrightarrow C^\infty_{p,\,q}(X,\,\C)$$

\noindent induced in every bidegree $(p,\,q)$ by any fixed Hermitian metric $\omega$ on $X$. (All the formal adjoints are computed w.r.t. the $L^2$ inner product defined by $\omega$ and so is the orthogonal projection $p'' = p''_\omega$ onto the harmonic space $\ker\Delta''$, where $\Delta'':=\bar\partial\bar\partial^\star + \bar\partial^\star\bar\partial$ is the usual $\bar\partial$-Laplacian induced by $\omega$.) This operator was introduced in [Pop16] where it was shown that every double class $[[\alpha^{p,\,q}]_{\bar\partial}]_{d_1}\in E_2^{p,\,q}(X)$ has a {\it unique} representative lying in the kernel of $\widetilde\Delta = \widetilde\Delta_\omega$ (cf. [Pop16, Theorem 1.1]).

\begin{Prop}\label{Prop:omega-injection} Let $X$ be a compact complex {\bf sGG manifold} with $\mbox{dim}_\C X=n$ and let $\omega$ be an arbitrary Hermitian metric on $X$. For any class $[[\alpha^{n-2,\,n}]_{\bar\partial}]_{d_1}\in E_2^{n-2,\,n}(X)$, let $\alpha^{n-2,\,n}_\omega$ be the {\bf $\widetilde\Delta_\omega$-harmonic} representative of $[[\alpha^{n-2,\,n}]_{\bar\partial}]_{d_1}$, let $\Omega^{n-1,\,n-1}_\omega\in C^\infty_{n-1,\,n-1}(X,\,\C)$ be the {\bf minimal $L^2_\omega$-norm} solution of the equation $\bar\partial\Omega^{n-1,\,n-1} = -\partial\alpha^{n-2,\,n}_\omega$ (cf. (\ref{eqn:choice-Omega})) and let $\beta^{n-2,\,n}_\omega\in C^\infty_{n-2,\,n}(X,\,\C)$ be the {\bf minimal $L^2_\omega$-norm} solution of the equation $\partial\beta^{n-2,\,n} = -\bar\partial\,\overline{\Omega^{n-1,\,n-1}_\omega}$ (cf. (\ref{eqn:choice-beta})).

The linear map $$j_\omega : E_2^{n-2,\,n}(X)\longrightarrow H^{2n-2}_{DR}(X,\,\C),  \hspace{3ex} j_\omega([[\alpha^{n-2,\,n}]_{\bar\partial}]_{d_1}) = \{\overline{\beta^{n-2,\,n}_\omega} + \Omega^{n-1,\,n-1}_\omega + \alpha^{n-2,\,n}_\omega\}_{DR},$$

\noindent is {\bf injective} and $T\circ j_\omega$ is the identity map of $E_2^{n-2,\,n}(X)$.

\end{Prop}

\noindent {\it Proof.} Suppose that $j_\omega([[\alpha^{n-2,\,n}]_{\bar\partial}]_{d_1}) = 0\in H^{2n-2}_{DR}(X,\,\C)$ for some $[[\alpha^{n-2,\,n}]_{\bar\partial}]_{d_1}\in E_2^{n-2,\,n}(X)$. Then, there exists a smooth $(2n-3)$-form $u=u^{n,\,n-3} + u^{n-1,\,n-2} + u^{n-2,\,n-1} + u^{n-3,\,n}$ such that $\overline{\beta^{n-2,\,n}_\omega} + \Omega^{n-1,\,n-1}_\omega + \alpha^{n-2,\,n}_\omega = du$. This implies that $\alpha^{n-2,\,n}_\omega = \partial u^{n-3,\,n} + \bar\partial u^{n-2,\,n-1}$. Since $\bar\partial u^{n-3,\,n}=0$, this further implies that $[[\alpha^{n-2,\,n}]_{\bar\partial}]_{d_1} = [[\alpha^{n-2,\,n}_\omega]_{\bar\partial}]_{d_1} = 0$. Thus, $j_\omega$ is injective.

The equality $T\circ j_\omega = \mbox{Id}_{E_2^{n-2,\,n}(X)}$ follows immediately from the definitions.  \hfill $\Box$

\subsection{Cohomology and sG metrics}\label{subsection:sG-metrics}

We now introduce strongly Gauduchon (sG) metrics into our discussion.

\begin{Def}\label{Def:cone_S} Let $X$ be a compact complex manifold with $\mbox{dim}_\C X=n$. 

(a)\, For every strongly Gauduchon metric (if any) $\omega>0$ on $X$, $\bar\partial\omega^{n-1}$ is $\partial$-exact. Let us denote by $\Gamma^{n-2,\,n}_\omega\in C^\infty_{n-2,\,n}(X,\,\C)$ the (unique) solution of minimal $L^2_\omega$-norm of the equation

\begin{equation}\label{eqn:sG-equation}\partial\Gamma^{n-2,\,n}_\omega = -\bar\partial\omega^{n-1}.\end{equation}

\noindent Since $\bar\partial\Gamma^{n-2,\,n}_\omega=0$ (for bidegree reasons) and $\partial\Gamma^{n-2,\,n}_\omega\in\mbox{Im}\,\bar\partial$, $\Gamma^{n-2,\,n}_\omega$ defines an element in $E_2^{n-2,\,n}(X)$.

 We consider the following subset

$${\cal S}_X:=\bigg\{\bigg[[\Gamma^{n-2,\,n}_\omega]_{\bar\partial}\bigg]_{d_1}\,\mid\, \omega \hspace{1ex} \mbox{is an sG metric on} \hspace{1ex} X\bigg\}\subset E_2^{n-2,\,n}(X)$$

\noindent that we call the {\bf $E_2sG$-cone} of $X$.

The real $(2n-2)$-form $\Gamma_\omega := \overline{\Gamma^{n-2,\,n}_\omega} + \omega^{n-1}  + \Gamma^{n-2,\,n}_\omega$ is $d$-closed, so it defines a real De Rham cohomology class $\{\Gamma_\omega\}_{DR}$. We consider the following subset

$$\widetilde{\cal S}_X:=\bigg\{\{\Gamma_\omega\}_{DR}\,\mid\, \omega \hspace{1ex} \mbox{is an sG metric on} \hspace{1ex} X\bigg\}\subset H^{2n-2}_{DR}(X,\,\R)$$

\noindent that we call the {\bf De Rham $E_2sG$-cone} of $X$.

\vspace{1ex}

(b)\, We also define the following variant of the $E_2sG$-cone of $X$ by dropping the $L^2_\omega$-norm minimality requirement on the solution $\Gamma^{n-2,\,n}$ of equation (\ref{eqn:sG-equation}): \begin{equation}\label{eqn:E2sG-cone-hat}\widehat{\cal S}_X = \bigg\{\bigg[[\Gamma^{n-2,\,n}]_{\bar\partial}\bigg]_{d_1}\,\mid\, \exists\,\omega \hspace{1ex} \mbox{Hermitian metric such that} \hspace{1ex} \partial\Gamma^{n-2,\,n} = -\bar\partial\omega^{n-1} \bigg\}\subset E_2^{n-2,\,n}(X).\end{equation} Every metric $\omega$ involved in (\ref{eqn:E2sG-cone-hat}) is sG, so we obviously have ${\cal S}_X\subset\widehat{\cal S}_X$. We do not know whether the reverse inclusion holds.

\end{Def}

The manifold $X$ is strongly Gauduchon if and only if ${\cal S}_X$ is non-empty.

\vspace{3ex}

We shall now prove that ${\cal S}_X$ (resp. $\widetilde{\cal S}_X$) is indeed a cone (i.e. a subset that is stable under multiplications by positive scalars) in $E_2^{n-2,\,n}(X)$ (resp. $H_{DR}^{2n-2}(X,\,\R)$). We need a few preliminaries.

\begin{Lem}\label{Lem:metric-rescaling}  Let $X$ be a compact complex manifold with $\mbox{dim}_\C X=n$. 

\vspace{1ex}

$(i)$\, For any Hermitian metric $\omega$ on $X$ and any positive real $\lambda$, the formal adjoints of $\bar\partial$ w.r.t. the metrics $\lambda\,\omega$ and $\omega$, as well as the corresponding $\bar\partial$-Laplacians, are related by the formulae

\begin{equation}\label{eqn:d-bar_star_relation}\bar\partial^\star_{\lambda\,\omega} = \frac{1}{\lambda}\,\bar\partial^\star_\omega \hspace{3ex} \mbox{and} \hspace{3ex} \Delta''_{\lambda\,\omega} = \frac{1}{\lambda}\,\Delta''_\omega   \end{equation}

\noindent in all bidegrees.

\vspace{1ex}

(ii)\, For any strongly Gauduchon metric $\omega$ on $X$ and any positive real $\lambda$, the forms $\Gamma^{n,\,n-2}_{\lambda\,\omega}:=\overline{\Gamma^{n-2,\,n}_{\lambda\,\omega}}$ and $\Gamma^{n,\,n-2}_{\omega}:=\overline{\Gamma^{n-2,\,n}_{\omega}}$ (see Definition \ref{Def:cone_S}) are related by the formula

$$\Gamma^{n,\,n-2}_{\lambda\,\omega} = \lambda^{n-1}\,\Gamma^{n,\,n-2}_{\omega}.$$

\noindent Consequently, we also have $\Gamma_{\lambda\,\omega} = \lambda^{n-1}\,\Gamma_{\omega}$ for any sG metric $\omega$ on $X$.

\end{Lem}

\noindent {\it Proof.} (i)\, Let us fix an arbitrary bidegree $(p,\,q)$. For any forms $\alpha, \beta$ of respective bidegrees $(p,\,q-1)$ and $(p,\,q)$, we have

$$\langle\langle\bar\partial\alpha,\,\beta\rangle\rangle_{\lambda\,\omega} = \lambda^n\,\int\limits_X\langle\bar\partial\alpha,\,\beta\rangle_{\lambda\,\omega}\,\frac{\omega^n}{n!} = \frac{\lambda^n}{\lambda^{p+q}}\,\int\limits_X\langle\bar\partial\alpha,\,\beta\rangle_{\omega}\,\frac{\omega^n}{n!} = \frac{\lambda^n}{\lambda^{p+q}}\,\langle\langle\bar\partial\alpha,\,\beta\rangle\rangle_{\omega} = \frac{\lambda^n}{\lambda^{p+q}}\,\langle\langle\alpha,\,\bar\partial^\star_\omega\beta\rangle\rangle_{\omega}$$

\noindent and

$$\langle\langle\alpha,\,\bar\partial^\star_{\lambda\,\omega}\beta\rangle\rangle_{\lambda\,\omega} = \frac{\lambda^n}{\lambda^{p+q-1}}\,\langle\langle\alpha,\,\bar\partial^\star_{\lambda\,\omega}\beta\rangle\rangle_{\omega}.$$

\noindent Since $\langle\langle\bar\partial\alpha,\,\beta\rangle\rangle_{\lambda\,\omega} = \langle\langle\alpha,\,\bar\partial^\star_{\lambda\,\omega}\beta\rangle\rangle_{\lambda\,\omega}$, the above formulae imply 

$$\frac{\lambda^n}{\lambda^{p+q}}\,\langle\langle\alpha,\,\bar\partial^\star_\omega\beta\rangle\rangle_{\omega} =\frac{\lambda^n}{\lambda^{p+q-1}}\,\langle\langle\alpha,\,\bar\partial^\star_{\lambda\,\omega}\beta\rangle\rangle_{\omega},  \hspace{3ex} \mbox{i.e.} \hspace{2ex} \langle\langle\alpha,\,\frac{1}{\lambda}\,\bar\partial^\star_\omega\beta\rangle\rangle_{\omega} = \langle\langle\alpha,\,\bar\partial^\star_{\lambda\,\omega}\beta\rangle\rangle_{\omega}$$

\noindent for all forms $\alpha$ and $\beta$. This proves the first formula in (\ref{eqn:d-bar_star_relation}).

 Since $\Delta'' = \bar\partial\bar\partial^\star + \bar\partial^\star\bar\partial$ (when $\Delta''$ and $\bar\partial^\star$ are computed w.r.t. the same metric), the latter formula in (\ref{eqn:d-bar_star_relation}) follows immediately from the former.  

\vspace{1ex}

(ii)\, By Definition \ref{Def:cone_S}, $\Gamma^{n,\,n-2}_{\omega}$ is the minimal $L^2_\omega$-norm solution of equation $\bar\partial\Gamma^{n,\,n-2} = -\partial\omega^{n-1}$. Consequently, the Neumann formula spells

\begin{equation}\label{eqn:Neumann}\Gamma^{n,\,n-2}_{\omega} = - \Delta_\omega^{''-1}\bar\partial^\star_\omega\,(\partial\omega^{n-1}). \end{equation}

\noindent Thus, we get: 

$$\Gamma^{n,\,n-2}_{\lambda\,\omega} = -\Delta_{\lambda\,\omega}^{''-1}\bar\partial^\star_{\lambda\,\omega}\,(\partial(\lambda\,\omega)^{n-1}) =-\lambda^{n-1}\,\Delta_{\omega}^{''-1}\bar\partial^\star_{\omega}\,(\partial\omega^{n-1}) = \lambda^{n-1}\,\Gamma^{n,\,n-2}_{\omega},$$

\noindent where we used the analogue of (\ref{eqn:Neumann}) for $\lambda\,\omega$ to get the first identity and (\ref{eqn:Neumann}) again to get the third identity.   \hfill $\Box$

\vspace{3ex}

\begin{Lem}\label{Lem:cone}  Let $X$ be a compact complex manifold with $\mbox{dim}_\C X=n$. The sets ${\cal S}_X$ and $\widehat{\cal S}_X$ are {\bf cones} in the $\C$-vector space $E_2^{n-2,\,n}(X)$, while the set $\widetilde{\cal S}_X$ is a cone in $H_{DR}^{2n-2}(X,\,\R)$. 

Moreover, the cone $\widehat{\cal S}_X$ is convex.

\end{Lem}

\noindent {\it Proof.} Let $[[\Gamma^{n-2,\,n}_\omega]_{\bar\partial}]_{d_1}\in{\cal S}_X$ and $\mu>0$ be arbitrary. Let $\lambda>0$ be the unique positive real such that $\lambda^{n-1}=\mu$. We have

$$\mu\,[[\Gamma^{n-2,\,n}_\omega]_{\bar\partial}]_{d_1} = [[\lambda^{n-1}\,\Gamma^{n-2,\,n}_\omega]_{\bar\partial}]_{d_1} = [[\Gamma^{n-2,\,n}_{\lambda\,\omega}]_{\bar\partial}]_{d_1},$$

\noindent where we used (ii) of Lemma \ref{Lem:metric-rescaling} to get the last identity. Now, $\lambda\,\omega$ is a strongly Gauduchon metric if $\omega$ is one, so $[[\Gamma^{n-2,\,n}_{\lambda\,\omega}]_{\bar\partial}]_{d_1}\in{\cal S}_X$. 

Consequently, ${\cal S}_X$ is stable under multiplications by positive scalars, hence it is a cone. The same goes for $\widetilde{\cal S}_X$ since (ii) of Lemma \ref{Lem:metric-rescaling} also applies to $\Gamma_\omega$. That $\widehat{\cal S}_X$ is a cone is trivial.

To prove the convexity of $\widehat{\cal S}_X$, it suffices to show that $\widehat{\cal S}_X$ is stable under additions. This is immediate since if $\partial\Gamma^{n-2,\,n}_i = -\bar\partial\omega_i^{n-1}$ for $i\in\{1,\,2\}$ and $\omega_i$ Hermitian metrics on $X$, then $\partial(\Gamma^{n-2,\,n}_1 + \Gamma^{n-2,\,n}_2) = -\bar\partial\omega_0^{n-1}$, where $\omega_0>0$ is the unique positive definite $C^\infty$ $(1,\,1)$-form on $X$ such that $\omega_0^{n-1}=\omega^{n-1}_1 + \omega^{n-1}_2>0$. Therefore, $[[\Gamma^{n-2,\,n}_1 + \Gamma^{n-2,\,n}_2]_{\bar\partial}]_{d_1}\in\widehat{\cal S}_X$.  \hfill $\Box$

\vspace{3ex}

The De Rham $E_2sG$-cone $\widetilde{\cal S}_X$ depends on the complex structure of $X$ and we now show this dependence to be lower semicontinuous in the sense described below in families of sGG manifolds. It actually suffices to assume that one fibre is sGG as all the nearby fibres are then sGG by the deformation openness of the sGG property (cf. [PU18, Corollary 1.7]).

\begin{Def}\label{Def:sG_H2-section} Let $\pi:{\cal X}\longrightarrow\Delta$ be a holomorphic family of compact complex $n$-dimensional manifolds over a ball $\Delta\subset\C^N$ centred at the origin. Suppose that the fibre $X_0:=\pi^{-1}(0)$ is {\bf strongly Gauduchon}. Let $X$ stand for the $C^\infty$ manifold that underlies the fibres $X_t$ with $t\in\Delta$. 

With every sG metric $\omega$ on $X_0$, we associate a {\bf local section $\tau_\omega$} of the constant real vector bundle ${\cal H}^{2n-2}_\R\longrightarrow\Delta$ whose fibre is the real De Rham cohomology space $H^{2n-2}_{DR}(X,\,\R)$ as follows.

As in Definition \ref{Def:cone_S}, we let $\Gamma_\omega := \overline{\Gamma^{n-2,\,n}_\omega} + \omega^{n-1}  + \Gamma^{n-2,\,n}_\omega$ be the real $d$-closed $(2n-2)$-form defined by the minimal $L^2_\omega$-norm solution $\Gamma^{n-2,\,n}_\omega$ of the equation $\partial\Gamma^{n-2,\,n}_\omega = -\bar\partial\omega^{n-1}$. (We put $\partial:=\partial_0$ and $\bar\partial:=\bar\partial_0$.) The component $(\Gamma_\omega)^{n-1,\,n-1}_t$ of $\Gamma_\omega$ of type $(n-1,\,n-1)$ for the complex structure of $X_t$ is positive definite if $t$ is close enough to $0$, by the continuity of the dependence on $t$ of $(\Gamma_\omega)^{n-1,\,n-1}_t$ and the positivity of $(\Gamma_\omega)^{n-1,\,n-1}_0 = \omega^{n-1}>0$. Hence, for every $t$ close to $0$, there exists a unique positive definite smooth $(1,\,1)$-form $\omega_t$ on $X_t$ such that $\omega_t^{n-1} = (\Gamma_\omega)^{n-1,\,n-1}_t >0$. In particular, $\omega_0 = \omega$.

Since $d\Gamma_\omega=0$, the form $\bar\partial_t\omega_t^{n-1}$ is $\partial_t$-exact (so $\omega_t$ is an sG metric on $X_t$). We let $\Gamma_{\omega_t}^{n-2,\,n}\in C^\infty_{n-2,\,n}(X_t,\,\C)$ be the minimal $L^2_{\omega_t}$-norm solution of the equation

\begin{equation}\label{eqn:Gamma_omega-t_n-2n-def}\partial_t\Gamma_{\omega_t}^{n-2,\,n} = -\bar\partial_t\omega_t^{n-1}\end{equation}

\noindent and we consider the real $d$-closed $(2n-2)$-form on $X$ defined as

$$\Gamma_\omega(t):=\overline{\Gamma_{\omega_t}^{n-2,\,n}} + \omega_t^{n-1} + \Gamma_{\omega_t}^{n-2,\,n}$$

\noindent for $t$ close to $0$. In particular, $\Gamma_\omega(0) = \Gamma_\omega$. Finally, we put

$$\tau_\omega(t):=\{\Gamma_\omega(t)\}_{DR}\in \widetilde{\cal S}_{X_t}\subset H^{2n-2}_{DR}(X,\,\R)$$

\noindent for all $t$ in a sufficiently small neighbourhood $U$ (depending on $\omega$) of $0$ in $\Delta$. 

\end{Def}

\vspace{2ex}

The {\bf lower semicontinuity} result for the De Rham $E_2sG$-cone $\widetilde{\cal S}_X\subset H^{2n-2}_{DR}(X,\,\R)$ when the complex structure of $X$ varies is the following

\begin{The}\label{The:E_2sG_semicontinuity} Let $\pi:{\cal X}\longrightarrow\Delta$ be a holomorphic family of compact complex $n$-dimensional manifolds over a ball $\Delta\subset\C^N$ centred at the origin. Suppose that the fibre $X_0:=\pi^{-1}(0)$ is an {\bf sGG manifold}. 

For every sG metric $\omega$ on $X_0$, the section $\tau_\omega$ of the constant real vector bundle ${\cal H}^{2n-2}_\R\longrightarrow\Delta$ on a small neighbourhood of $0$ in $\Delta$ constructed in Definition \ref{Def:sG_H2-section} is $C^\infty$.

\end{The}

\vspace{2ex}

In particular, every element $\{\Gamma_\omega\}_{DR}$ of the De Rham $E_2sG$-cone $\widetilde{\cal S}_{X_0}$ of $X_0$ extends to a $C^\infty$ family of elements $\{\Gamma_\omega(t)\}_{DR}$ of the De Rham $E_2sG$-cones $\widetilde{\cal S}_{X_t}$ of the nearby fibres $X_t$. Moreover, there is such an extension for every representative $\Gamma_\omega$ of the given De Rham class $\{\Gamma_\omega\}_{DR}$ defined by an sG metric $\omega$. So, in this sense, the De Rham $E_2sG$-cone $\widetilde{\cal S}_{X_0}$ of $X_0$ can only be ``smaller'' than the De Rham $E_2sG$-cones $\widetilde{\cal S}_{X_t}$ of the nearby fibres $X_t$.

\vspace{2ex}

\noindent {\it Proof.} By the well-known Neumann formula, the minimal $L^2_{\omega_t}$-norm solution of equation (\ref{eqn:Gamma_omega-t_n-2n-def}) is

$$\Gamma_{\omega_t}^{n-2,\,n} = -(\partial_t)^\star_{\omega_t}\Delta_{\omega_t}^{'-1}(\bar\partial_t\omega_t^{n-1}),$$

\noindent where $(\partial_t)^\star_{\omega_t}$ is the formal adjoint of $\partial_t$ w.r.t. the $L^2$ inner product induced by the metric $\omega_t$, while $\Delta'_{\omega_t}=\partial_t(\partial_t)^\star_{\omega_t} + (\partial_t)^\star_{\omega_t}\partial_t$ is the $\partial$-Laplacian induced by $\omega_t$ and $\Delta_{\omega_t}^{'-1}$ stands for its Green operator.

Now, the $(n-1,\,n)$-form $\bar\partial_t\omega_t^{n-1}$ varies in a $C^\infty$ way with $t$ and so do the differential operators $\Delta'_{\omega_t}$ and $(\partial_t)^\star_{\omega_t}$. Moreover, the classical Kodaira-Spencer theory (cf. [KS60]) applied to the $C^\infty$ family $(\Delta'_{\omega_t})_{t\in\Delta}$ of {\it elliptic} differential operators acting in bidegree $(n-1,\,n)$ ensures that the family 
$(\Delta_{\omega_t}^{'-1})_{t\in\Delta}$ of their Green operators is again $C^\infty$ if the dimensions of the kernels $\ker\Delta'_{\omega_t}$ (that are isomorphic to the $\partial$-cohomology spaces $H^{n-1,\,n}_{\partial}(X_t,\,\C)$ by the Hodge isomorphism) are independent of $t$. However, $H^{n-1,\,n}_{\partial}(X_t,\,\C)$ is $\C$-anti-linearly isomorphic to $H^{n,\,n-1}_{\bar\partial}(X_t,\,\C)$ by conjugation, while the latter vector space is Serre-dual to $H^{0,\,1}_{\bar\partial}(X_t,\,\C)$, so its dimension equals the Hodge number $h^{0,\,1}_{\bar\partial}(t)$ of the fibre $X_t$ for every $t$.

Here is where the sGG assumption on the fibre $X_0$ comes in. By [PU18, Corollary 1.7], it ensures that the Hodge numbers $h^{0,\,1}_{\bar\partial}(t)$ are independent of $t$ when $t$ varies in a small enough neighbourhood of $0$. Thus, the Green operators $\Delta_{\omega_t}^{'-1}$ in bidegree $(n-1,\,n)$, hence also the $(n-2,\,n)$-forms $\Gamma_{\omega_t}^{n-2,\,n}$, vary in a $C^\infty$ way with $t$ near $0$. Since so also do (for trivial reasons) the $(n-1,\,n-1)$-forms $\omega_t^{n-1}$, we infer that the smooth $(2n-2)$-forms 

$$\Gamma_\omega(t):=\overline{\Gamma_{\omega_t}^{n-2,\,n}} + \omega_t^{n-1} + \Gamma_{\omega_t}^{n-2,\,n}$$

\noindent vary in a $C^\infty$ way with $t$ in a small enough neighbourhood of $0$. The application of the De Rham cohomology class being a smooth operation, we conclude that $\tau_\omega(t):=\{\Gamma_\omega(t)\}_{DR}$ depends in a $C^\infty$ way on $t$ varying in a small enough neighbourhood of $0\in\Delta$.  \hfill $\Box$

\subsection{Real-valued cohomology}\label{subsection:R-valued}

We shall now deal with the real version of some of the objects introduced in $\S.$\ref{subsection:C-valued} for the sake of enhanced flexibility.

\begin{Def}\label{Def:potential_realE2} Let $X$ be a compact complex manifold with $\mbox{dim}_\C X=n$.

\vspace{1ex}

(a)\, For any element $[[\alpha^{n-2,\,n}]_{\bar\partial}]_{d_1}\in E_2^{n-2,\,n}(X)$ and any representative $\alpha^{n-2,\,n}\in C^\infty_{n-2,\,n}(X,\,\C)$ of this double class, we know from (\ref{eqn:E_2_rep}) that $\bar\partial\alpha=0$ (trivial here for bidegree reasons) and that there exists a (not necessarily real and non-unique) form $\Omega^{n-1,\,n-1}\in C^\infty_{n-1,\,n-1}(X,\,\C)$ such that $\partial\alpha^{n-2,\,n} = -\bar\partial\Omega^{n-1,\,n-1}$. 

We refer to any such form $\Omega^{n-1,\,n-1}$ as an {\bf $(n-1,\,n-1)$-potential} of $\alpha^{n-2,\,n}$.

\vspace{1ex}

(b)\, We define the {\bf real} part $E_2^{n-2,\,n}(X)_\R$ of the $\C$-vector space $E_2^{n-2,\,n}(X)$ by selecting the classes representable by forms admitting a {\bf real} $(n-1,\,n-1)$-potential:

$$E_2^{n-2,\,n}(X)_\R:=\bigg\{[[\alpha^{n-2,\,n}]_{\bar\partial}]_{d_1}\in E_2^{n-2,\,n}(X)\,\mid\,\exists\alpha^{n-2,\,n} \hspace{1ex} \mbox{representative having a real potential}\hspace{1ex}\Omega^{n-1,\,n-1} \bigg\}.$$

\noindent By definition, $E_2^{n-2,\,n}(X)_\R$ is a {\bf real} vector subspace of $E_2^{n-2,\,n}(X)$. 

\end{Def}

We shall now consider the real version of the map $T$ introduced in $\S.$\ref{subsection:C-valued}.

\begin{Lem}\label{Lem:T_R} Let $X$ be a compact complex manifold with $\mbox{dim}_\C X=n$.

\vspace{1ex}

(i)\, The following inclusion holds: $E_2^{n-2,\,n}(X)_\R\subset\ker d_2^{n-2,\,n}$.

\vspace{1ex}

(ii)\, The restriction to $H_{DR}^{2n-2}(X,\,\R)$ of the map $T$ defined in Proposition \ref{Prop:T_def-image}, namely the map

$$T_\R: H^{2n-2}_{DR}(X,\,\R)\longrightarrow E_2^{n-2,\,n}(X)_\R, \hspace{3ex}  \{\alpha\}_{DR}\mapsto\bigg[[\alpha^{n-2,\,n}]_{\bar\partial}\bigg]_{d_1},$$

\noindent assumes its values in the {\bf real} space $E_2^{n-2,\,n}(X)_\R$ and is {\bf surjective}.

\end{Lem}

\noindent {\it Proof.} (i)\, Let $[[\alpha^{n-2,\,n}]_{\bar\partial}]_{d_1}\in E_2^{n-2,\,n}(X)_\R$ with $\partial\alpha^{n-2,\,n} = -\bar\partial\Omega^{n-1,\,n-1}$ for some {\bf real} $(n-1,\,n-1)$-form $\Omega^{n-1,\,n-1}$. Thanks to (\ref{eqn:E_2_map_def}), we have

$$d_2([[\alpha^{n-2,\,n}]_{\bar\partial}]_{d_1}) = -[[\partial\Omega^{n-1,\,n-1}]_{\bar\partial}]_{d_1} = [[\bar\partial\overline{\alpha^{n-2,\,n}}]_{\bar\partial}]_{d_1} =0\in E_2^{n,\,n-1}(X),$$

\noindent because by conjugating the identity defining $\Omega^{n-1,\,n-1}$ and using the fact that $\Omega^{n-1,\,n-1}$ is real, we get $\partial\Omega^{n-1,\,n-1} = -\bar\partial\overline{\alpha^{n-2,\,n}}$. 

Therefore, $[[\alpha^{n-2,\,n}]_{\bar\partial}]_{d_1}\in\ker d_2^{n-2,\,n}$.

\vspace{1ex}

(ii)\, Let $\{\alpha\}_{DR}\in H_{DR}^{2n-2}(X,\,\R)$ and pick a real representative $\alpha = \alpha^{n,\,n-2} +  \alpha^{n-1,\,n-1} +  \alpha^{n-2,\,n}$. Since $\alpha$ is real, $\alpha^{n-1,\,n-1}$ is real. Since $\alpha$ is $d$-closed, $\partial\alpha^{n-2,\,n} = -\bar\partial\alpha^{n-1,\,n-1}$. Thus, $\alpha^{n-1,\,n-1}$ is a real $(n-1,\,n-1)$-potential of $\alpha^{n-2,\,n}$, so $T(\{\alpha\}_{DR}) = [[\alpha^{n-2,\,n}]_{\bar\partial}]_{d_1}\in E_2^{n-2,\,n}(X)_\R$.

To prove that $T$ is surjective, let $[[\alpha^{n-2,\,n}]_{\bar\partial}]_{d_1}\in E_2^{n-2,\,n}(X)_\R$ with $\partial\alpha^{n-2,\,n} = -\bar\partial\Omega^{n-1,\,n-1}$ for some {\bf real} $(n-1,\,n-1)$-form $\Omega^{n-1,\,n-1}$. Then, the $(2n-2)$-form

$$\alpha:=\overline{\alpha^{n-2,\,n}} + \Omega^{n-1,\,n-1} + \alpha^{n-2,\,n}$$

\noindent is real, $d$-closed and $T(\{\alpha\}_{DR}) = [[\alpha^{n-2,\,n}]_{\bar\partial}]_{d_1}\in E_2^{n-2,\,n}(X)_\R$.  \hfill $\Box$

\vspace{2ex}

We are now in a position to see that the real space $E_2^{n-2,\,n}(X)_\R$ contains the $E_2sG$-cone $\widehat{\cal S}_X$ of $X$ defined in $\S.$\ref{subsection:C-valued} as an open cone.

\begin{Lem}\label{Lem:open-cone} Let $X$ be a compact complex manifold with $\mbox{dim}_\C X=n$. The inclusions ${\cal S}_X\subset\widehat{\cal S}_X\subset E_2^{n-2,\,n}(X)_\R$ hold and the cone $\widehat{\cal S}_X$ is {\bf open} in $E_2^{n-2,\,n}(X)_\R$.

\end{Lem}

\noindent {\it Proof.} The inclusion ${\cal S}_X\subset\widehat{\cal S}_X$ is obvious and has already been noticed. Let $[[\Gamma^{n-2,\,n}]_{\bar\partial}]_{d_1}\in\widehat{\cal S}_X$ be arbitrary. So, there exists a representative $\Gamma^{n-2,\,n}\in C^\infty_{n-2,\,n}(X,\,\C)$ of this $E_2$-class and an sG metric $\omega$ on $X$ such that $\partial\Gamma^{n-2,\,n} = -\bar\partial\omega^{n-1}$. Since $\omega^{n-1}$ is {\bf real}, $[[\Gamma^{n-2,\,n}]_{\bar\partial}]_{d_1}\in E_2^{n-2,\,n}(X)_\R$. This proves the inclusion $\widehat{\cal S}_X\subset E_2^{n-2,\,n}(X)_\R$.

Let $[[\alpha^{n-2,\,n}]_{\bar\partial}]_{d_1}\in E_2^{n-2,\,n}(X)_\R$ and let $\varepsilon>0$. Then, there exists a {\bf real} form $\Omega^{n-1,\,n-1}\in C^\infty_{n-1,\,n-1}(X,\,\C)$ such that $\bar\partial\Omega^{n-1,\,n-1} = -\partial\alpha^{n-2,\,n}$. We get

$$\partial(\Gamma^{n-2,\,n} + \varepsilon\,\alpha^{n-2,\,n}) = -\bar\partial(\omega^{n-1} + \varepsilon\,\Omega^{n-1,\,n-1}).$$

\noindent On the other hand, the $(n-1,\,n-1)$-form $\omega^{n-1} + \varepsilon\,\Omega^{n-1,\,n-1}$ is real for every $\varepsilon$ and is positive definite if $\varepsilon>0$ is small enough. Therefore, for every small $\varepsilon>0$, there exists a unique positive definite $(1,\,1)$-form $\rho_\varepsilon>0$ such that $\rho_\varepsilon^{n-1} =\omega^{n-1} + \varepsilon\,\Omega^{n-1,\,n-1}$. We get $\partial\rho_\varepsilon^{n-1}= - \bar\partial(\overline{\Gamma^{n-2,\,n}} + \varepsilon\,\overline{\alpha^{n-2,\,n}})$, hence $\partial\rho_\varepsilon^{n-1}$ is $\bar\partial$-exact, so $\rho_\varepsilon$ is a strongly Gauduchon metric on $X$. Consequently, 

$$[[\Gamma^{n-2,\,n}]_{\bar\partial}]_{d_1} + \varepsilon\,[[\alpha^{n-2,\,n}]_{\bar\partial}]_{d_1}\in\widehat{\cal S}_X$$

\noindent for every small $\varepsilon>0$. This proves that $\widehat{\cal S}_X$ is open in $E_2^{n-2,\,n}(X)_\R$.  \hfill $\Box$

\subsection{Duality of positive cones in the $E_2$ cohomology}\label{subsection:cone-duality}

It was proved in [PU18b] (by means of the pseudo-differential Laplacian $\widetilde\Delta$ introduced in [Pop16] that gives a Hodge theory for the second page of the Fr\"olicher spectral sequence) that for any compact complex $n$-dimensional manifold $X$ and every $p,q\in\{0,\dots , n\}$, the canonical bilinear pairing 

\begin{equation}\label{eqn:duality_E2}E_2^{p,\,q}(X)\times E_2^{n-p,\,n-q}(X)\longrightarrow\C, \hspace{6ex} \bigg([[\alpha]_{\bar\partial}]_{d_1},\,[[\beta]_{\bar\partial}]_{d_1}\bigg)\mapsto \int\limits_X\alpha\wedge\beta,\end{equation}

\noindent is well defined (i.e. independent of the choices of representatives of the $E_2$-cohomology classes involved) and non-degenerate. Hence, it defines a Serre-type {\bf duality} between $E_2^{p,\,q}(X)$ and $E_2^{n-p,\,n-q}(X)$.

Under this duality, the closure of our $E_2sG$ cone $\widehat{\cal S}_X\subset E_2^{n-2,\,n}(X)$, consisting of those $E_2$-classes of type $(n-2,\,n)$ that are ``positive'' in the sense of Definition \ref{Def:cone_S}, has a dual cone in $E_2^{2,\,0}(X)$ that we will now describe. To this end, we will introduce ad hoc notions of {\it real} and {\it positive} $(2,\,0)$-forms and currents that run counter to the standard definitions of real and positive forms and currents of bidegree $(p,\,p)$, but propose a not so far-fetched analogue thereof in this bidegree that is relevant to holomorphic symplectic geometry. A possible extension of this geometry on sGG manifolds is one of our motivations and will hopefully be attempted in future work. 

 In this subsection, we will establish an $E_2$ analogue for the bidegrees $(2,\,0)$ and $(n-2,\,n)$ of Lamari's duality (cf. [Lam99, lemme 3.3]) between Demailly's pseudo-effective cone ${\cal E}(X)\subset H^{1,\,1}_{BC}(X,\,\R)$ (consisting of the Bott-Chern cohomology classes of all the $d$-closed, positive $(1,\,1)$-currents $T\geq 0$ on $X$, see [Dem92]) and the closure of the Gauduchon cone ${\cal G}_X\subset H^{n-1,\,n-1}_A(X,\,\R)$ introduced in [Pop15] (consisting of the Aeppli cohomology classes of all the Gauduchon metrics $\omega^{n-1}>0$ on $X$).

We will assume throughout this subsection that $X$ is an sGG manifold. This will guarantee that every Gauduchon metric is actually strongly Gauduchon (see [PU18a]).

\begin{Def}\label{Def:positivity_skewed-bideg} Let $X$ be an {\bf sGG} compact complex $n$-dimensional manifold. 

\vspace{1ex}

(i)\, We consider the following sets: \begin{eqnarray}\nonumber V & = & \bigg\{\Gamma^{n-2,\,n}\in C^\infty_{n-2,\,n}(X,\,\C)\,\mid\, \exists\,\omega \hspace{1ex} \mbox{Hermitian metric such that} \hspace{1ex} \partial\Gamma^{n-2,\,n} = -\bar\partial\omega^{n-1} \bigg\}, \\
\nonumber E & = & \bigg\{\Gamma^{n-2,\,n}\in C^\infty_{n-2,\,n}(X,\,\C)\,\mid\, \partial\Gamma^{n-2,\,n}\in\mbox{Im}\,\bar\partial\bigg\}, \\
\nonumber E_{\R} & = & \bigg\{\Gamma^{n-2,\,n}\in C^\infty_{n-2,\,n}(X,\,\C)\,\mid\,\exists\, \Omega^{n-1,\,n-1} \hspace{1ex} \mbox{{\bf real} form such that} \hspace{1ex} \partial\Gamma^{n-2,\,n} = -\bar\partial\Omega^{n-1,\,n-1}\bigg\}.\end{eqnarray}

\vspace{1ex}

\noindent Thus, $V\subset E_{\R} \subset E\subset C^\infty_{n-2,\,n}(X,\,\C)$ and $E$ consists of the smooth $(n-2,\,n)$-forms that are $E_2$-closed (their $\bar\partial$-closedness is automatic for bidegree reasons), while $E_{\R}$ consists of the {\bf real} (in this ad hoc sense) such forms and $V$ consists of the {\bf positive} (in this ad hoc sense) such forms. Note that any metric $\omega$ featuring in the definition of $V$ is automatically strongly Gauduchon (or, equivalently, Gauduchon since $X$ is assumed sGG). 

\vspace{1ex}

(ii)\, Fix an arbitrary Hermitian metric $\gamma$ on $X$. Let $p_{Im\,\bar\partial}^{(\gamma)} : C^\infty_{n-2,\,n}(X,\,\C)\longrightarrow \mbox{Im}\,\bar\partial$ be the orthogonal projection w.r.t. the $L^2_\gamma$-inner product onto the closed subspace of $\bar\partial$-exact $(n-2,\,n)$-forms, induced by the standard Hodge-theoretical $L^2_\gamma$-orthogonal $3$-space decomposition

$$C^\infty_{n-2,\,n}(X,\,\C) = \ker\Delta''\oplus\mbox{Im}\,\bar\partial\oplus\mbox{Im}\,\bar\partial^\star.$$

\vspace{1ex}

We consider the following sets: \begin{eqnarray}\nonumber U_\gamma & = & \bigg\{\Gamma^{n-2,\,n}\in C^\infty_{n-2,\,n}(X,\,\C)\,\mid\, \exists\,\omega \hspace{1ex} \mbox{Hermitian metric such that} \hspace{1ex} p_{Im\,\bar\partial}^{(\gamma)}(\partial\Gamma^{n-2,\,n}) = -\bar\partial\omega^{n-1} \bigg\}, \end{eqnarray}

\noindent $C^\infty_{n-2,\,n}(X,\,\R)_{\gamma}  = $ \\

$ \bigg\{\alpha^{n-2,\,n}\in C^\infty_{n-2,\,n}(X,\,\C)\,\mid\, \exists\,\beta^{n-1,\,n-1} \hspace{1ex} \mbox{{\bf real} form such that} \hspace{1ex} p_{Im\,\bar\partial}^{(\gamma)}(\partial\alpha^{n-2,\,n}) = -\bar\partial\beta^{n-1,\,n-1}\bigg\}.$

\vspace{3ex}

\noindent Thus, $U_\gamma\subset C^\infty_{n-2,\,n}(X,\,\R)_{\gamma}\subset C^\infty_{n-2,\,n}(X,\,\C)$ and $U_\gamma$ consists of the {\bf $\gamma$-positive} (in this ad hoc sense) smooth $(n-2,\,n)$-forms, while $C^\infty_{n-2,\,n}(X,\,\R)_{\gamma}$ consists of the {\bf $\gamma$-real} (in this ad hoc sense) such forms. Unlike the sets defined under (i), these sets are not subjected to any $E_2$-closedness condition. In particular, the following inclusions hold:

$$V\subset U_\gamma \hspace{3ex} \mbox{and} \hspace{3ex} E_\R\subset C^\infty_{n-2,\,n}(X,\,\R)_{\gamma}$$

\noindent for every Hermitian metric $\gamma$ on $X$.

\vspace{1ex}

(iii)\, In the context of (ii), by a {\bf $\gamma$-real current of bidegree $(2,\,0)$} on $X$ we mean any continuous $\R$-linear form

$$\tau^{2,\,0} : C^\infty_{n-2,\,n}(X,\,\R)_{\gamma} \longrightarrow\R.$$

\noindent By such a current being {\bf $\gamma$-positive} we mean that $\tau^{2,\,0}$ evaluates non-negatively on every element of $U_\gamma$. (So, in particular, the zero current $\tau^{2,\,0}= 0$ is $\gamma$-positive.)

By a $\gamma$-real current $\tau^{2,\,0}$ of bidegree $(2,\,0)$ being {\bf $E_2$-exact} we mean that $\tau^{2,\,0}$ vanishes identically on the $\R$-vector space $E_\R$ of ``real'' $E_2$-closed $(n-2,\,n)$-forms defined under (i). \footnote{This last notion is in keeping with the usual duality according to which a current is exact (w.r.t. a given cohomology) if and only if it vanishes identically on the closed (w.r.t the same cohomology) $C^\infty$ forms of complementary bidegree.}

\end{Def}

The following properties of the above sets are immediate to check.

\begin{Lem}\label{Lem:sets_properties} (a)\, The set $E$ is a {\bf closed} $\C$-vector subspace of $C^\infty_{n-2,\,n}(X,\,\C)$, the sets $E_\R$ and $C^\infty_{n-2,\,n}(X,\,\R)_{\gamma}$ are {\bf closed} $\R$-vector subspaces of $C^\infty_{n-2,\,n}(X,\,\C)$, while $V$ and $U_\gamma$ are {\bf open convex} cones in $E_\R$, respectively $C^\infty_{n-2,\,n}(X,\,\R)_{\gamma}$.

\vspace{1ex}

(b)\, The following identity holds: $$U_\gamma\cap E_\R = V.$$

\end{Lem}

\noindent {\it Proof.} (a)\, The closedness conclusion for $E$ follows from the well-known fact (itself a consequence of standard elliptic theory on compact manifolds) that $\mbox{Im}\,\bar\partial$ is {\it closed} in the space of $C^\infty$ forms in which it lies. 

To see that $E_\R$ is closed in $C^\infty_{n-2,\,n}(X,\,\C)$, we need one further step. Let $\Gamma_j^{n-2,\,n}\rightarrow\Gamma^{n-2,\,n}\in C^\infty_{n-2,\,n}(X,\,\C)$ in the $C^\infty$ topology as $j\rightarrow +\infty$, where $\Gamma_j^{n-2,\,n}\in E_\R$ for every $j\in\N$. For every $j$, let $\Gamma_j^{n-1,\,n-1} = -\Delta^{''-1}_\gamma\bar\partial^\star_\gamma\partial\Gamma_j^{n-2,\,n}\in C^\infty_{n-1,\,n-1}(X,\,\C)$ be the solution of minimal $L^2_\gamma$-norm of the equation $\bar\partial\Gamma_j^{n-1,\,n-1} = -\partial\Gamma_j^{n-2,\,n}$. (So, the set of all the solutions is the affine subspace $\Gamma_j^{n-1,\,n-1} + \ker\bar\partial\subset C^\infty_{n-1,\,n-1}(X,\,\C)$ and the hypothesis $\Gamma_j^{n-2,\,n}\in E_\R$ means that $(\Gamma_j^{n-1,\,n-1} + \ker\bar\partial)\cap C^\infty_{n-1,\,n-1}(X,\,\R)\neq\emptyset$, where $C^\infty_{n-1,\,n-1}(X,\,\R)\subset C^\infty_{n-1,\,n-1}(X,\,\C)$ is the real vector subspace of {\it real} forms.) Then $\mbox{Im}\,\bar\partial\ni\partial\Gamma_j^{n-2,\,n}\rightarrow\partial\Gamma^{n-2,\,n}$ in the $C^\infty$ topology as $j\rightarrow +\infty$, so $\partial\Gamma^{n-2,\,n}\in\mbox{Im}\,\bar\partial$ since $\mbox{Im}\,\bar\partial$ is {\it closed}. Moreover, $\Gamma^{n-1,\,n-1} = -\Delta^{''-1}_\gamma\bar\partial^\star_\gamma\partial\Gamma^{n-2,\,n}\in C^\infty_{n-1,\,n-1}(X,\,\C)$ is the solution of minimal $L^2_\gamma$-norm of the equation $\bar\partial\Gamma^{n-1,\,n-1} = -\partial\Gamma^{n-2,\,n}$, so $\Gamma_j^{n-1,\,n-1}\rightarrow\Gamma^{n-1,\,n-1}$ in the $C^\infty$ topology as $j\rightarrow +\infty$ because the restriction to $\mbox{Im}\,\bar\partial$ of the operator $\Delta^{''-1}_\gamma\bar\partial^\star_\gamma$ is continuous in the $C^\infty$ topology. Since $(\Gamma_j^{n-1,\,n-1} + \ker\bar\partial)\cap C^\infty_{n-1,\,n-1}(X,\,\R)\neq\emptyset$ for all $j\in\N$ and $C^\infty_{n-1,\,n-1}(X,\,\R)$ is {\it closed} in $C^\infty_{n-1,\,n-1}(X,\,\C)$, we get $(\Gamma^{n-1,\,n-1} + \ker\bar\partial)\cap C^\infty_{n-1,\,n-1}(X,\,\R)\neq\emptyset$. This means that $\Gamma^{n-2,\,n}\in E_\R$. Thus, $E_\R$ is closed in $C^\infty_{n-2,\,n}(X,\,\C)$.

The closedness of $C^\infty_{n-2,\,n}(X,\,\R)_{\gamma}$ in $C^\infty_{n-2,\,n}(X,\,\C)$ can be proved in the same way since the projection $p_{Im\,\bar\partial}^{(\gamma)}$ is continuous w.r.t. the $C^\infty$ topology.

The convexity of $V$ and $U_\gamma$ follows from the linearity of the operators $\partial$, $\bar\partial$ and $p_{Im\,\bar\partial}^{(\gamma)}$ involved in their definitions and from the convexity of the set of Gauduchon metrics (itself a consequence of the existence of a unique positive definite $(n-1)^{st}$ root for every positive definite $(n-1,\,n-1)$-form). 

Let us prove that $U_\gamma$ is open in $C^\infty_{n-2,\,n}(X,\,\R)_{\gamma}$. (The openness of $V$ in $E_\R$ can be proved in a similar way.) Let $\Gamma^{n-2,\,n}\in U_\gamma$ and $\alpha^{n-2,\,n}\in C^\infty_{n-2,\,n}(X,\,\R)_{\gamma}$ be arbitrary. By definition, there exist a Hermitian metric $\omega$ and a real form $\beta^{n-1,\,n-1}\in C^\infty_{n-1,\,n-1}(X,\,\R)$ such that

$$p_{Im\,\bar\partial}^{(\gamma)}(\partial\Gamma^{n-2,\,n}) = -\bar\partial\omega^{n-1}   \hspace{3ex} \mbox{and} \hspace{3ex} p_{Im\,\bar\partial}^{(\gamma)}(\partial\alpha^{n-2,\,n}) = -\bar\partial\beta^{n-1,\,n-1}.$$

\noindent Thus, for every constant $\varepsilon>0$, we get $p_{Im\,\bar\partial}^{(\gamma)}(\partial(\Gamma^{n-2,\,n} + \varepsilon\,\alpha^{n-2,\,n})) = -\bar\partial(\omega^{n-1} + \varepsilon\,\beta^{n-1,\,n-1}).$ Since $\beta^{n-1,\,n-1}$ is real and $\omega^{n-1}$ is positive definite, $\omega^{n-1} + \varepsilon\,\beta^{n-1,\,n-1}$ is positive definite for all sufficiently small $\varepsilon>0$. Therefore, $\Gamma^{n-2,\,n} + \varepsilon\,\alpha^{n-2,\,n}\in U_\gamma$. This proves that $U_\gamma$ is open in $C^\infty_{n-2,\,n}(X,\,\R)_{\gamma}$.  

\vspace{1ex}

(b)\, To prove the inclusion ``$\subset$'', let $\Gamma^{n-2,\,n}\in U_\gamma\cap E_\R$. Since $\Gamma^{n-2,\,n}\in E_\R$, $\partial\Gamma^{n-2,\,n}\in\mbox{Im}\,\bar\partial$, so $p_{Im\,\bar\partial}^{(\gamma)}(\partial\Gamma^{n-2,\,n}) = \partial\Gamma^{n-2,\,n}$. Thus, $\partial\Gamma^{n-2,\,n} =-\bar\partial\omega^{n-1}$ for some Hermitian metric $\omega$ thanks to $\Gamma^{n-2,\,n}$ lying in $U_\gamma$. Therefore, $\Gamma^{n-2,\,n}\in V$. This proves the inclusion ``$\subset$''. The reverse inclusion is obvious.  \hfill $\Box$

\vspace{2ex} 

The last preliminary remark that we make serves as an example pointing out a very particular way of constructing {\it real-valued} linear maps on $E_\R$. Such maps occur below in a more general form. (See the last hypothesis of Proposition \ref{Prop:cone-duality}.)  

\begin{Lem}\label{Lem:real-valued_necessary-cond} If a form $\theta^{2,\,0}\in C^\infty_{2,\,0}(X,\,\C)$ is of the shape $\theta^{2,\,0}=\partial\xi^{1,\,0}$ such that the $(1,\,1)$-form $\bar\partial\xi^{1,\,0}$ is {\bf real}, then $\int_X\theta^{2,\,0}\wedge\Gamma^{n-2,\,n}$ is {\bf real} for every $\Gamma^{n-2,\,n}\in E_\R$.

\end{Lem}

\noindent {\it Proof.} Let $\Gamma^{n-2,\,n}\in E_\R$. Then, $\partial\Gamma^{n-2,\,n} = -\bar\partial\Omega^{n-1,\,n-1}$ for some {\it real} $(n-1,\,n-1)$-form $\Omega^{n-1,\,n-1}$. Applying Stokes's Theorem twice, we get

$$\int\limits_X\theta^{2,\,0}\wedge\Gamma^{n-2,\,n} = \int\limits_X\partial\xi^{1,\,0}\wedge\Gamma^{n-2,\,n} = - \int\limits_X\xi^{1,\,0}\wedge\bar\partial\Omega^{n-1,\,n-1} = - \int\limits_X\bar\partial\xi^{1,\,0}\wedge\Omega^{n-1,\,n-1}.$$

\noindent The last quantity is real since both forms $\bar\partial\xi^{1,\,0}$ and $\Omega^{n-1,\,n-1}$ are real.  \hfill $\Box$

\vspace{2ex}

We are now in a position to prove the duality result we have been aiming for. Both the statement and the proof parallel those of Lemma 3.3 in [Lam99].

\begin{Prop}\label{Prop:cone-duality}  Let $X$ be an {\bf sGG} compact complex $n$-dimensional manifold on which an arbitrary Hermitian metric $\gamma$ has been fixed. Let $\theta^{2,\,0}\in C^\infty_{2,\,0}(X,\,\C)$ satisfy the condition

$$\int\limits_X\theta^{2,\,0}\wedge\Gamma^{n-2,\,n}\geq 0$$

\noindent for every $\Gamma^{n-2,\,n}\in C^\infty_{n-2,\,n}(X,\,\C)$ for which there exists a Hermitian metric $\omega$ on $X$ such that $\partial\Gamma^{n-2,\,n} = -\bar\partial\omega^{n-1}$ (i.e. for every $\Gamma^{n-2,\,n}\in V$). Suppose, moreover, that $\int_X\theta^{2,\,0}\wedge\Gamma^{n-2,\,n}\in\R$ for every $\Gamma^{n-2,\,n}\in E_\R$. 

Then, there exists a $\gamma$-positive current $\tau^{2,\,0} : C^\infty_{n-2,\,n}(X,\,\R)_{\gamma}\longrightarrow\R$ of bidegree $(2,\,0)$ on $X$ such that $\theta^{2,\,0} - \tau^{2,\,0}$ is $E_2$-exact in the sense that it vanishes identically on $E_\R$.

\end{Prop}

\vspace{2ex}

Note that if we assume $\theta^{2,\,0}$ to be $E_2$-closed (i.e. $\bar\partial\theta^{2,\,0}=0$ and $\partial\theta^{2,\,0}\in\mbox{Im}\,\bar\partial$, which in bidegree $(2,\,0)$ is equivalent to assuming that $d\theta^{2,\,0}=0$), it defines a class $[[\theta^{2,\,0}]_{\bar\partial}]_{d_1}\in E_2^{2,\,0}$ and the integral $\int_X\theta^{2,\,0}\wedge\Gamma^{n-2,\,n}$ is independent of the choice of representative of this class. Thus, Proposition \ref{Prop:cone-duality} implies that the dual of the closure of the $E_2sG$-cone $\widehat{\cal S}_X$ defined in (\ref{eqn:E2sG-cone-hat}) under the duality $E_2^{2,\,0}(X)\times E_2^{n-2,\,n}(X)\longrightarrow\C$ is the closed convex cone in $E_2^{2,\,0}(X)$ consisting of the $E_2$-classes $[[\theta^{2,\,0}]_{\bar\partial}]_{d_1}$ ``representable'' by $\gamma$-positive currents $\tau^{2,\,0} : C^\infty_{n-2,\,n}(X,\,\R)_{\gamma}\longrightarrow\R$.

\vspace{3ex}

\noindent {\it Proof of Proposition \ref{Prop:cone-duality}.} We follow closely Lamari's arguments of the proof of Lemma 3.3. in [Lam99]. The form $\theta^{2,\,0}$ defines a $\C$-linear map 

$$\theta^{2,\,0}:C^\infty_{n-2,\,n}(X,\,\C)\longrightarrow\C,  \hspace{3ex} \Gamma^{n-2,\,n}\mapsto\int\limits_X\theta^{2,\,0}\wedge\Gamma^{n-2,\,n}.$$ 

\noindent The hypothesis imposed on $\theta^{2,\,0}$ translates to $\theta^{2,\,0}_{|V}\geq 0$. Thus, there are two cases. 

\vspace{2ex}

{\it Case $1$.} Suppose there exists $\Gamma^{n-2,\,n}_0\in V\subset E_\R$ such that $\int_X\theta^{2,\,0}\wedge\Gamma^{n-2,\,n}_0= 0$. This implies that $\theta^{2,\,0}_{|E_\R}\equiv 0$.

Indeed, fix an arbitrary $\Gamma^{n-2,\,n}\in E_\R$ and let $\Gamma^{n-2,\,n}_t:=(1-t)\,\Gamma^{n-2,\,n}_0 + t\,\Gamma^{n-2,\,n}$ for $t\in\R$. Then $\Gamma^{n-2,\,n}_t\in E_\R$ for all $t\in[0,\,1]$ since $E_\R$ is convex. Moreover,  for all $t\in\R$, we have

$$f(t):=\int\limits_X\theta^{2,\,0}\wedge\Gamma^{n-2,\,n}_t = (1-t)\,\int\limits_X\theta^{2,\,0}\wedge\Gamma^{n-2,\,n}_0 + t\,\int\limits_X\theta^{2,\,0}\wedge\Gamma^{n-2,\,n} = t\,\int\limits_X\theta^{2,\,0}\wedge\Gamma^{n-2,\,n}.$$

\noindent In particular, $f(0)=0$. Meanwhile, $V$ is open in $E_\R$ (cf. Lemma \ref{Lem:sets_properties}) and $\Gamma^{n-2,\,n}_0\in V$, so $\Gamma^{n-2,\,n}_t\in V$ for all $t$ close enough to $0$. Since $\theta^{2,\,0}_{|V}\geq 0$, we infer that $f(t)\geq 0$ for all $t\in[-\varepsilon,\,\varepsilon]$ for some small $\varepsilon>0$. This means that $t\,\int_X\theta^{2,\,0}\wedge\Gamma^{n-2,\,n}\geq 0$ for all $t\in[-\varepsilon,\,\varepsilon]$, which is impossible unless $\int_X\theta^{2,\,0}\wedge\Gamma^{n-2,\,n}=0$. This proves that $\theta^{2,\,0}_{|E_\R}\equiv 0$, so we can choose $\tau^{2,\,0}=0$ (which is $\gamma$-positive).

\vspace{2ex}

{\it Case $2$.} Suppose that $\theta^{2,\,0}_{|V}> 0$. Let $F\subset E_\R$ be the kernel of the restriction $\theta^{2,\,0}_{|E_\R} : E_\R\longrightarrow \R$. Thus, $F$ has real codimension $1$ in $E_\R$ and

$$U_\gamma\cap F=\emptyset.$$

\noindent To see the last identity, suppose there exists $\Gamma^{n-2,\,n}\in U_\gamma\cap F$. Then, $\int_X\theta^{2,\,0}\wedge\Gamma^{n-2,\,n}=0$ because $\Gamma^{n-2,\,n}\in F=\ker(\theta^{2,\,0}_{|E_\R})$. Meanwhile, $\int_X\theta^{2,\,0}\wedge\Gamma^{n-2,\,n}>0$ because $F\subset E_\R$, so $\Gamma^{n-2,\,n}\in U_\gamma\cap E_\R = V$ (see (b) of Lemma \ref{Lem:sets_properties} for the last identity) and $\theta^{2,\,0}_{|V}> 0$. This is a contradiction.

Since $U_\gamma$ is a convex open subset of $C^\infty_{n-2,\,n}(X,\,\R)_\gamma$, $F$ is a convex closed subset of $C^\infty_{n-2,\,n}(X,\,\R)_\gamma$ and $U_\gamma$ and $F$ are disjoint, the Hahn-Banach Separation Theorem allows us to separate them. Consequently, there exists a continuous $\R$-linear form

$$l^{2,\,0}:C^\infty_{n-2,\,n}(X,\,\R)_\gamma\longrightarrow\R$$

\noindent such that $l^{2,\,0}_{|U_\gamma}>0$ and  $l^{2,\,0}_{|F}=0$.

The first condition implies that the $\gamma$-real current $l^{2,\,0}$ of bidegree $(2,\,0)$ is $\gamma$-positive. 

Let $\Gamma_1^{n-2,\,n}\in V$. Then $\int_X\theta^{2,\,0}\wedge\Gamma^{n-2,\,n}_1>0$ and $\int_X l^{2,\,0}\wedge\Gamma^{n-2,\,n}_1>0$, so there exists a constant $\lambda>0$ such that $\int_X\theta^{2,\,0}\wedge\Gamma^{n-2,\,n}_1 = \lambda\,\int_X l^{2,\,0}\wedge\Gamma^{n-2,\,n}_1$. This means that $(\theta^{2,\,0} - \lambda\,l^{2,\,0})_{|\R\,\Gamma^{n-2,\,n}_1}=0$. But we also have $(\theta^{2,\,0} - \lambda\,l^{2,\,0})_{|F}=0$. Since the real codimension of $F$ in $E_\R$ is $1$, we get that $(\theta^{2,\,0} - \lambda\,l^{2,\,0})_{|E_\R}=0$. If we put $\tau^{2,\,0}:=\lambda\,l^{2,\,0}$, we are done.  \hfill $\Box$

\section{The $h$-$\partial\bar\partial$ property of compact complex manifolds}\label{section:h-dd-bar}

Let $X$ be a compact complex manifold with $\mbox{dim}_\C X=n$. We now consider the {\it adiabatic limit} construction of the differential operator $d_h=h\partial + \bar\partial$ (cf. (\ref{eqn:d_h_def})) that was introduced in [Pop17] for every constant $h>0$. Allowing now $h$ to be negative, some obvious properties include the following: \begin{eqnarray}\label{eqn:obvious-prop}\nonumber & (i)& \overline{d}_h = h\,d_{h^{-1}};  \hspace{3ex}  (ii)\, \overline{d}_{-h} = -h\,d_{-h^{-1}};\\
\nonumber &(iii)& d_{h_1}d_{h_2} = (h_1-h_2)\,\partial\bar\partial;  \hspace{5ex} \mbox{in particular,}   \hspace{2ex} d_h\,d_{-h^{-1}} = (h+ \frac{1}{h})\,\partial\bar\partial;\\
    &(iv)& \frac{h+1}{h^2+1}\,d_h + \frac{h(h-1)}{h^2+1}\,d_{-h^{-1}} =d, \end{eqnarray} 

\noindent for all $h\in\R\setminus\{0\}$.

When a Hermitian metric $\omega$ has been fixed on $X$, the formal adjoint $d_h^\star$ of $d_h$ w.r.t. $\omega$ induces together with $d_h$ a Laplace-type operator in the usual way: $$\Delta_h:d_h\,d_h^\star + d_h^\star d_h:C^\infty_k(X,\,\C)\longrightarrow C^\infty_k(X,\,\C),$$ \noindent for every $k\in\{0,\dots , 2n\}$. This {\bf h-Laplacian} is elliptic (cf. [Pop17]). Identity (ii) in (\ref{eqn:obvious-prop}) implies \begin{eqnarray}\label{eqn:obvious-prop_bis}\overline{\Delta}_{-h} = h^2\,\Delta_{-h^{-1}},   \hspace{3ex} \mbox{for all} \hspace{1ex} h\in\R\setminus\{0\}.\end{eqnarray} We shall now continue the study of the operators $d_h$, both from a metric and an intrinsic angle.

\subsection{Commutation relations and BKN identity for the operators $d_h$}\label{subsection:d_h-commutation}

Let us fix an arbitrary Hermitian metric $\omega$ on $X$. All the formal adjoints will be computed w.r.t. $\omega$, as will the (pointwise and formal) adjoint $\Lambda = \Lambda_\omega$ of the multiplication operator $\L = \L_\omega :=\omega\wedge\cdot$. Recall the standard torsion operator of type $(1,\,0)$ (cf. [Dem84]) $\tau =[\Lambda,\,\partial\omega\wedge\cdot]$ and the Hermitian commutation relations (cf. again [Dem84]):

$$\partial^\star + \tau^\star = i\,[\Lambda,\,\bar\partial]  \hspace{3ex} \mbox{and}  \hspace{3ex} \bar\partial^\star + \bar\tau^\star = -i\,[\Lambda,\,\partial].$$

We will infer the following

\begin{Lem}\label{Lem:h-commutation-rel} Let $(X,\,\omega)$ be a complex Hermitian manifold. For every $h\in\R\setminus\{0\}$, we define the {\bf $h$-torsion operator} of type $(1,\,0)$ induced by $\omega$ by $\tau_h:=[\Lambda,\,d_h\omega\wedge\cdot]$.

The following {\bf Hermitian $h$-commutation relations} hold on differential forms of any degree: \begin{eqnarray}\nonumber (a)\,(d_h+\tau_h)^\star = -i\,[\Lambda,\, \overline{d}_{-h}];  \hspace{3ex} & &  \hspace{3ex}  (b)\,(\overline{d}_h+\overline{\tau}_h)^\star = i\,[\Lambda,\, d_{-h}];\\
\nonumber (c)\,d_h+\tau_h = i\,[\overline{d}_{-h}^\star,\,\omega\wedge\cdot];  \hspace{3ex} & &  \hspace{3ex}  (d)\,\overline{d}_h+\overline{\tau}_h = -i\,[d_{-h}^\star,\,\omega\wedge\cdot].\end{eqnarray}

\end{Lem} 

\noindent {\it Proof.} Since (b) is the conjugate of (a), while the implications $(a)\implies (c)$ and $(b) \implies (d)$ are obtained by taking adjoints, it suffices to prove (a).

Using the above definitions and the standard Hermitian commutation relations, we get \begin{eqnarray}\nonumber d_h^\star & = & h\partial^\star + \bar\partial^\star = i\,[\Lambda,\,h\bar\partial] -h\tau^\star -i\,[\Lambda,\,\partial] -\bar\tau^\star = -i\,[\Lambda,\,\overline{d}_{-h}] - (h\tau^\star + \bar\tau^\star)\end{eqnarray} \noindent and \begin{eqnarray}\nonumber h\tau^\star + \bar\tau^\star & = & [(h\,\partial\omega\wedge\cdot)^\star,\,\omega\wedge\cdot] + [(\bar\partial\omega\wedge\cdot)^\star,\,\omega\wedge\cdot] = [(d_h\omega\wedge\cdot)^\star,\,\omega\wedge\cdot] = \tau_h^\star.\end{eqnarray}  \noindent Summing up these identities, we get (a).  \hfill $\Box$

\vspace{2ex}

An immediate consequence is the following

\begin{Cor}\label{Cor:rough-BKN} Let $(X,\,\omega)$ be a complex Hermitian manifold. For every $h\in\R\setminus\{0\}$, the following {\bf rough $h$-Bochner-Kodaira-Nakano (h-BKN) identity} holds on differential forms of any degree:

$$\Delta_h = \overline{\Delta}_{-h} + [\overline{d}_{-h},\,\overline{\tau}_{-h}^\star] - [d_h,\,\tau_h^\star].$$

\end{Cor}

\noindent {\it Proof.} Using the $h$-commutation relation (a) of Lemma \ref{Lem:h-commutation-rel} for the second identity below, we get

$$\Delta_h=[d_h,\,d_h^\star] = -i\,[d_h,\,[\Lambda,\,\overline{d}_{-h}]] - [d_h,\,\tau_h^\star].$$ 

 On the other hand, the Jacobi identity spells:

$$-[d_h,\,[\Lambda,\,\overline{d}_{-h}]] + [\Lambda,\,[\overline{d}_{-h},\,d_h]] + [\overline{d}_{-h},\,[d_h,\,\Lambda]] =0.$$

\noindent Since $[\overline{d}_{-h},\,d_h]=0$ whenever $h\neq 0$, the second term above vanishes. Meanwhile, $[d_h,\,\Lambda] = i\,(\overline{d}_{-h} + \overline{\tau}_{-h})^\star$ as follows from the $h$-commutation relation (b) of Lemma \ref{Lem:h-commutation-rel} after replacing $h$ with $-h$. Therefore, we get $-i\,[d_h,\,[\Lambda,\,\overline{d}_{-h}]] = [\overline{d}_{-h},\,(\overline{d}_{-h} + \overline{\tau}_{-h})^\star] = \overline{\Delta}_{-h} + [\overline{d}_{-h},\,\overline{\tau}_{-h}^\star]$ and the formula follows.  \hfill $\Box$

\vspace{3ex}

Another immediate consequence is the following anti-commutation statement in the K\"ahler case.

\begin{Cor}\label{Cor:Kaehler_anti-comm} Let $(X,\,\omega)$ be a compact {\bf K\"ahler} manifold. For every $h\in\R\setminus\{0\}$, the following identities hold:

$$[d_h,\,d^\star_{-h^{-1}}] = 0  \hspace{3ex} \mbox{and} \hspace{3ex} [d_{-h^{-1}},\, d_h^\star] = 0.$$

\end{Cor}

\noindent {\it Proof.} The latter identity is the adjoint of the former, so it suffices to prove the former one. When $h$ has been replaced by $-h^{-1}$, the $h$-commutation relation (a) of Lemma \ref{Lem:h-commutation-rel} spells $d^\star_{-h^{-1}} = -i\,[\Lambda,\,\overline{d}_{h^{-1}}] = -\frac{i}{h}\,[\Lambda,\,d_h] $ since $\tau_h=0$ for every $h$ when $\omega$ is K\"ahler and identity (i) in (\ref{eqn:obvious-prop}) has been used to infer the last identity. Therefore, $[d_h,\,d^\star_{-h^{-1}}] = -\frac{i}{h}\,[d_h,\,[\Lambda,\,d_h]]$ when $\omega$ is K\"ahler.

Now, the Jacobi identity yields:

$$-[d_h,\,[\Lambda,\,d_h]] + [\Lambda,\,[d_h,\,d_h]] + [d_h,\,[d_h,\,\Lambda]] =0.$$

\noindent Since $[d_h,\,d_h] =0$ and $[d_h,\,\Lambda] = -[\Lambda,\,d_h]$, we get $[d_h,\,[\Lambda,\,d_h]]=0$. 

Consequently, $[d_h,\,d^\star_{-h^{-1}}] = -\frac{i}{h}\,[d_h,\,[\Lambda,\,d_h]] =0$ and we are done.  \hfill $\Box$

\vspace{3ex}

Taking our cue from [Dem84], we shall now refine the above BKN formula by incorporating the $1^{st}$ order terms into a twisted Laplace-type operator on the r.h.s. of the identity so that the discrepancy terms become of order zero. We begin with some preliminary computations.

\begin{Lem}\label{Lem:preliminary_BKN-refined} Let $(X,\,\omega)$ be a complex Hermitian manifold. For every $h\in\R\setminus\{0\}$, the following identities hold: \\

(i)\, $[L,\,\tau_h] = 3\,d_h\omega\wedge\cdot$,  \hspace{6ex}   (ii)\, $[\Lambda,\,\tau_h] = 2i\,\overline{\tau}_{-h}^\star$,  \hspace{6ex}    (iii)\, $[d_h,\,\overline{d}_{-h}^\star] = -[d_h,\,\overline{\tau}_{-h}^\star]$, \\

 (iv)\, $[d_h,\,d_h^\star] + [d_h,\,\tau_h^\star] -[\overline{d}_{-h},\,\overline{\tau}_{-h}^\star]  = [d_h + \tau_h,\,d_h^\star + \tau_h^\star] + S_\omega^{(h)},$ where

$$S_\omega^{(h)} : = \frac{i}{2}\,[\Lambda,\,[\Lambda,\,\overline{d}_{-h}d_h\omega\wedge\cdot]] - [d_h\omega\wedge\cdot,\, (d_h\omega\wedge\cdot)^\star].$$

\end{Lem}

\noindent {\it Proof.} (i)\, The definition of $\tau_h$ and the Jacobi identity yield the first and respectively the second identities below: 

$$[L,\,\tau_h] = [L,\,[\Lambda,\,d_h\omega\wedge\cdot]] = - [\Lambda,\,[d_h\omega\wedge\cdot,\,L]] - [d_h\omega\wedge\cdot,\,[L,\,\Lambda]].$$ 

\noindent Now, $[d_h\omega\wedge\cdot,\,L] = d_h\omega\wedge(\omega\wedge\cdot) - \omega\wedge d_h\omega\wedge\cdot =0$, so the first term on the r.h.s. above vanishes. Meanwhile, it is standard that $[L,\,\Lambda]=(k-n)\,\mbox{Id}$ on $k$-forms. So for any $k$-form $u$, we get 

$$[d_h\omega\wedge\cdot,\,[L,\,\Lambda]]\,u = d_h\omega\wedge([L,\,\Lambda]\,u) - [L,\,\Lambda]\,(d_h\omega\wedge u) = (k-n)\,d_h\omega\wedge u - (k+3-n)\,d_h\omega\wedge u = -3\,d_h\omega\wedge u.$$ 

\noindent Thus, $[d_h\omega\wedge\cdot,\,[L,\,\Lambda]] = -3\,d_h\omega\wedge\cdot$ and (i) follows.

\vspace{1ex}

(ii)\, We know from the $h$-commutation relation (c) of Lemma \ref{Lem:h-commutation-rel} that $\tau_h = i\,[\overline{d}_{-h}^\star,\,\omega\wedge\cdot] - d_h$. Hence, using also (b) of Lemma \ref{Lem:h-commutation-rel}, we get

$$[\Lambda,\,\tau_h] = i\,[\Lambda,\,[\overline{d}_{-h}^\star,\,\omega\wedge\cdot]] - [\Lambda,\,d_h] = i\,[\Lambda,\,[\overline{d}_{-h}^\star,\,\omega\wedge\cdot]] + i\,(\overline{d}_{-h} + \overline{\tau}_{-h})^\star.$$

\noindent The Jacobi identity spells

$$[\Lambda,\,[\overline{d}_{-h}^\star,\,\omega\wedge\cdot]] + [\overline{d}_{-h}^\star,\,[\omega\wedge\cdot,\,\Lambda]] + [\omega\wedge\cdot,\,[\Lambda,\,\overline{d}_{-h}^\star]] =0.$$

\noindent Since $[\omega\wedge\cdot,\,\Lambda] = (k-n)\,\mbox{Id}$ on $k$-forms, we get $[\overline{d}_{-h}^\star,\,[\omega\wedge\cdot,\,\Lambda]] = \overline{d}_{-h}^\star$. Meanwhile, $[\omega\wedge\cdot,\,[\Lambda,\,\overline{d}_{-h}^\star]] = [[\overline{d}_{-h},\,\omega\wedge\cdot],\,\Lambda]^\star$, so we get

$$[\Lambda,\,\tau_h] = -i\,[[\overline{d}_{-h},\,\omega\wedge\cdot],\,\Lambda]^\star - i\,\overline{d}_{-h}^\star + i\,(\overline{d}_{-h} + \overline{\tau}_{-h})^\star.$$

Moreover, for an arbitrary form $u$, we get \begin{eqnarray}\nonumber[\overline{d}_{-h},\,\omega\wedge\cdot]\,u  = (\partial -h\bar\partial)\,(\omega\wedge u) - \omega\wedge(\partial u -h\bar\partial u) = (\partial\omega - h\,\bar\partial\omega)\wedge u = \overline{d}_{-h}\omega\wedge u.\end{eqnarray}

\noindent Thus, $[\overline{d}_{-h},\,\omega\wedge\cdot] =  \overline{d}_{-h}\,\omega\wedge\cdot$ and we finally get \begin{eqnarray}\nonumber[\Lambda,\,\tau_h] & = & -i\,[\overline{d}_{-h}\,\omega\wedge\cdot,\,\Lambda]^\star - i\,\overline{d}_{-h}^\star + i\,(\overline{d}_{-h} + \overline{\tau}_{-h})^\star \\
\nonumber & = & i\,\overline{\tau}_{-h}^\star - i\,\overline{d}_{-h}^\star + i\,(\overline{d}_{-h} + \overline{\tau}_{-h})^\star = 2i\,\overline{\tau}_{-h}^\star,\end{eqnarray}

\noindent where the second identity followed from the definition of $\tau_h$ by replacing $h$ with $-h$ and then taking conjugates and adjoints.

This proves (ii).

\vspace{1ex}

(iii)\, The Jacobi identity yields

$$-[d_h,\,[\Lambda,\,d_h]] + [\Lambda,\,[d_h,\,d_h]] + [d_h,\,[d_h,\,\Lambda]] =0.$$

\noindent Since $[d_h,\,d_h]=0$ (because $d_h^2=0$), and $[d_h,\,\Lambda] = -[\Lambda,\,d_h]$, we get $[d_h,\,[\Lambda,\,d_h]]=0$. Using the $h$-commutation relation (b) of Lemma \ref{Lem:h-commutation-rel}, this means that $[d_h,\,\overline{\tau}_{-h}^\star] = - [d_h,\,\overline{d}_{-h}^\star]$, or equivalently that $[\tau_h,\,\overline{d}_{-h}^\star] = -[d_h,\,\overline{d}_{-h}^\star],$ where the latter identity was obtained from the former by taking adjoints, conjugates and replacing $h$ with $-h$. In other words, we have \begin{equation}\label{eqn:A3}[d_h,\,\overline{d}_{-h}^\star] = -[\tau_h,\,\overline{d}_{-h}^\star] = -[d_h,\,\overline{\tau}_{-h}^\star].\end{equation}

This proves (iii). 

\vspace{1ex}

(iv)\, Applying part (ii) and then the Jacobi identity, we get \begin{eqnarray}\label{eqn:BKN_preliminary0}[\overline{d}_{-h},\,\overline{\tau}_{-h}^\star] = -\frac{i}{2}\,[\overline{d}_{-h},\,[\Lambda,\,\tau_h]] = -\frac{i}{2}\,[\Lambda,\,[\tau_h,\,\overline{d}_{-h}]] -\frac{i}{2}\,[\tau_h,\,[\overline{d}_{-h},\,\Lambda]].\end{eqnarray}

\noindent On the other hand, \begin{eqnarray}\nonumber[\tau_h,\,\overline{d}_{-h}] & \stackrel{(a)}{=} & [\overline{d}_{-h},\,\tau_h] \stackrel{(b)}{=} [\overline{d}_{-h},\,[\Lambda,\,d_h\omega\wedge\cdot]] \stackrel{(c)}{=} [\Lambda,\,[d_h\omega\wedge\cdot,\,\overline{d}_{-h}]] + [d_h\omega\wedge\cdot,\,[\overline{d}_{-h},\, \Lambda]] \\ 
\nonumber & \stackrel{(d)}{=} & [\Lambda,\,\overline{d}_{-h}d_h\omega\wedge\cdot] - i\,[d_h\omega\wedge\cdot,\, d_h^\star + \tau_h^\star],\end{eqnarray}

\noindent where (a) follows from $\tau_h$ and $\overline{d}_{-h}$ being operators of odd degrees, (b) follows from the definition of $\tau_h$, (c) follows from the Jacobi identity, while the latter term in (d) follows from the $h$-commutation relation (b) of Lemma \ref{Lem:h-commutation-rel} and the former term in (d) follows from the following easy computation:

$$[d_h\omega\wedge\cdot,\,\overline{d}_{-h}]\,u = d_h\omega\wedge \overline{d}_{-h}u + \overline{d}_{-h}(d_h\omega\wedge u) = \overline{d}_{-h}d_h\omega\wedge u,$$

\noindent for any form $u$. 

Taking the bracket with $\Lambda$ in the above formula for $[\tau_h,\,\overline{d}_{-h}]$, we get

\begin{eqnarray}\label{eqn:BKN_preliminary1}[\Lambda,\,[\tau_h,\,\overline{d}_{-h}]] = [\Lambda,\,[\Lambda,\,\overline{d}_{-h}d_h\omega\wedge\cdot]] -i\, [\Lambda,\, [d_h\omega\wedge\cdot,\, d_h^\star + \tau_h^\star]].\end{eqnarray}

\noindent Applying again the Jacobi formula for the last term, we get \begin{eqnarray}\label{eqn:BKN_preliminary2}\nonumber[\Lambda,\, [d_h\omega\wedge\cdot,\, d_h^\star + \tau_h^\star]] & = & - [d_h\omega\wedge\cdot,\,[d_h^\star + \tau_h^\star,\,\Lambda]] + [d_h^\star + \tau_h^\star,\,[\Lambda,\,d_h\omega\wedge\cdot]] \\
\nonumber & = & - [d_h\omega\wedge\cdot,\, [\omega\wedge\cdot,\,d_h+\tau_h]^\star] + [d_h^\star + \tau_h^\star,\,\tau_h],\\
 & = & -2\,[d_h\omega\wedge\cdot,\, (d_h\omega\wedge\cdot)^\star] + [d_h^\star + \tau_h^\star,\,\tau_h],   \end{eqnarray}

\noindent where the first term on the last line is given by the following simple computation. For any form $u$, we have $[\omega\wedge\cdot,\,d_h]\,u = \omega\wedge d_h u - d_h(\omega\wedge u) = -d_h\omega\wedge u$. Thus, $[\omega\wedge\cdot,\,d_h] = -d_h\omega\wedge\cdot$. Combined with identity (i), this yields $[\omega\wedge\cdot,\,d_h+\tau_h] = 2\,d_h\omega\wedge\cdot$.

Putting (\ref{eqn:BKN_preliminary1}) and (\ref{eqn:BKN_preliminary2}) together, we get

$$[\Lambda,\,[\tau_h,\,\overline{d}_{-h}]] = [\Lambda,\,[\Lambda,\,\overline{d}_{-h}d_h\omega\wedge\cdot]] + 2i\,[d_h\omega\wedge\cdot,\, (d_h\omega\wedge\cdot)^\star] -i\,[d_h^\star + \tau_h^\star,\,\tau_h],$$

\noindent which, in turn, combines with (\ref{eqn:BKN_preliminary0}) to yield

$$[\overline{d}_{-h},\,\overline{\tau}_{-h}^\star]  =  -\frac{i}{2}\,[\Lambda,\,[\Lambda,\,\overline{d}_{-h}d_h\omega\wedge\cdot]] + [d_h\omega\wedge\cdot,\, (d_h\omega\wedge\cdot)^\star] -\frac{1}{2}\,[d_h^\star + \tau_h^\star,\,\tau_h] -\frac{i}{2}\,[\tau_h,\,[\overline{d}_{-h},\,\Lambda]].$$

\noindent Since $-i\,[\Lambda,\,\overline{d}_{-h}] = d_h^\star + \tau_h^\star$ by the $h$-commutation relation (a) of Lemma \ref{Lem:h-commutation-rel}, we get

$$-[\overline{d}_{-h},\,\overline{\tau}_{-h}^\star]  = [d_h^\star + \tau_h^\star,\,\tau_h] + \frac{i}{2}\,[\Lambda,\,[\Lambda,\,\overline{d}_{-h}d_h\omega\wedge\cdot]] - [d_h\omega\wedge\cdot,\, (d_h\omega\wedge\cdot)^\star].$$

\noindent Adding $[d_h,\,d_h^\star] + [d_h,\,\tau_h^\star]$ on either side of the above identity, we get

$$[d_h,\,d_h^\star] + [d_h,\,\tau_h^\star] -[\overline{d}_{-h},\,\overline{\tau}_{-h}^\star]  = [d_h + \tau_h,\,d_h^\star + \tau_h^\star] + S_\omega^{(h)},$$

\noindent where $S_\omega^{(h)} : = \frac{i}{2}\,[\Lambda,\,[\Lambda,\,\overline{d}_{-h}d_h\omega\wedge\cdot]] - [d_h\omega\wedge\cdot,\, (d_h\omega\wedge\cdot)^\star]$. This proves (iv).  \hfill $\Box$

\vspace{3ex}

We can now state the main result of this subsection.

\begin{The}\label{Thm:BKN} Let $(X,\,\omega)$ be a complex Hermitian manifold. For every $h\in\R\setminus\{0\}$, the following {\bf refined $h$-Bochner-Kodaira-Nakano (h-BKN) identity} holds on differential forms of any degree: 

$$\Delta_h = [\overline{d}_{-h} + \overline{\tau}_{-h},\,\overline{d}_{-h}^\star + \overline{\tau}_{-h}^\star] + T_\omega^{(h)},$$

\noindent where $T_\omega^{(h)}$ is the zero-th order operator defined by $$T_\omega^{(h)}:= -\frac{i}{2}\,[\Lambda,\,[\Lambda,\,d_h\overline{d}_{-h}\omega\wedge\cdot]] - [\overline{d}_{-h}\omega\wedge\cdot,\, (\overline{d}_{-h}\omega\wedge\cdot)^\star].$$

In particular, if the metric $\omega$ is {\bf K\"ahler}, $d_h\omega=0$ hence $\tau_h=0$ and $T_\omega^{(h)}=0$, so we get $$\Delta_h = \overline{\Delta}_{-h} \hspace{2ex} (\mbox{for every} \hspace{1ex} h\in\R) \hspace{3ex} \mbox{and}  \hspace{3ex} \Delta_h = h^2\,\Delta_{-h^{-1}}  \hspace{2ex} (\mbox{for every} \hspace{1ex} h\in\R\setminus\{0\}).$$ 

\noindent The latter identity follows from the former thanks to (\ref{eqn:obvious-prop_bis}).

\end{The}

\noindent {\it Proof.} Combining (iv) of Lemma \ref{Lem:preliminary_BKN-refined} with the rough BKN formula of Corollary \ref{Cor:rough-BKN}, we get

$$\Delta_h + [d_h + \tau_h,\,d_h^\star + \tau_h^\star] + S_\omega^{(h)} = \overline{\Delta}_{-h} + [\overline{d}_{-h},\,\overline{\tau}_{-h}^\star] - [d_h,\,\tau_h^\star] + [d_h,\,d_h^\star] + [d_h,\,\tau_h^\star] -[\overline{d}_{-h},\,\overline{\tau}_{-h}^\star].$$

\noindent Since $[d_h,\,d_h^\star] = \Delta_h$, the last formula reduces to

$$\overline{\Delta}_{-h} = [d_h + \tau_h,\,d_h^\star + \tau_h^\star] + S_\omega^{(h)}.$$

The refined h-BKN identity follows from this by taking conjugates and replacing $h$ with $-h$. \hfill $\Box$

\subsection{$h$-$\partial\bar\partial$-manifolds}\label{subsection:h-dd-bar}

The standard $\partial\bar\partial$-lemma asserts that every compact K\"ahler manifold is a $\partial\bar\partial$-manifold. We will now investigate the analogue of this statement in our $d_h$-cohomology context.

\begin{The}\label{Thm:Laplace_sum} Let $(X,\,\omega)$ be a compact {\bf K\"ahler} manifold. As usual, we let $\Delta:=dd^\star + d^\star d$. For every $h\in\R\setminus\{0\}$, the following identity holds:

$$\Delta = \frac{(h+1)^2}{(h^2+1)^2}\,\Delta_h + \frac{(h-1)^2}{(h^2+1)^2}\,(h^2\,\Delta_{-h^{-1}}).$$

\end{The}

\noindent {\it Proof.} Using (iv) of (\ref{eqn:obvious-prop}) and the obvious identity $\Delta_h=  [d_h,\,d_h^\star]$ for every $h$, we get \begin{eqnarray}\nonumber\Delta  = [d,\,d^\star] & = & \frac{(h+1)^2}{(h^2+1)^2}\,\Delta_h + \frac{h^2\,(h-1)^2}{(h^2+1)^2}\,\Delta_{-h^{-1}} \\
\nonumber & + & \frac{(h+1)h(h-1)}{(h^2+1)^2}\,[d_h,\,d^\star_{-h^{-1}}] + \frac{(h+1)h(h-1)}{(h^2+1)^2}\,[d_{-h^{-1}},\,d^\star_h].\end{eqnarray}

\noindent Since the metric $\omega$ is supposed to be K\"ahler, $[d_h,\,d^\star_{-h^{-1}}]=0$ and $[d_{-h^{-1}},\,d^\star_h]=0$ by Corollary \ref{Cor:Kaehler_anti-comm}. The statement follows.  \hfill $\Box$

\vspace{2ex}

An immediate consequence of Theorems \ref{Thm:BKN} and \ref{Thm:Laplace_sum} is the following proportionality statement.

\begin{Cor}\label{Cor:Laplace_proportionality} Let $(X,\,\omega)$ be a compact {\bf K\"ahler} manifold. For every $h\in\R\setminus\{0\}$, the following identities hold on differential forms of any degree:

$$\Delta = \frac{2}{h^2+1}\,\Delta_h = \frac{2h^2}{h^2+1}\,\Delta_{-h^{-1}} = \frac{2}{h^2+1}\,\overline{\Delta}_{-h}.$$

\end{Cor}

\vspace{3ex}

We pause briefly to notice that the above proportionality statement reproves, in conjunction with the main result of [Pop17], the standard fact that the K\"ahler property of compact complex manifolds implies the degeneration at $E_1$ of the Fr\"olicher spectral sequence. Yet another proof will be implicit further down by putting together Theorems \ref{The:h-dd-bar} and \ref{The:h-dd-bar_E1}.

\begin{Cor}\label{Cor:Kaehler_E_1_inf} (standard) Let $(X,\,\omega)$ be a compact {\bf K\"ahler} manifold. The Fr\"olicher spectral sequence of $X$ degenerates at $E_1$.

\end{Cor}

\noindent {\it Proof.} We know from Corollary \ref{Cor:Laplace_proportionality} that $\Delta_h =\frac{h^2+1}{2}\,\Delta$ for every $h\in\R\setminus\{0\}$ in every degree $k$. In particular, $\ker\Delta_h = \ker\Delta$ for all $h\neq 0$. Let $\delta_h^{(k)}>0$ be the smallest positive eigenvalue of $\Delta_h:C^\infty_k(X,\,\C)\longrightarrow C^\infty_k(X,\,\C)$ acting on $k$-forms and let $u_h^{(k)}\in C^\infty_k(X,\,\C)$ be a corresponding eigenvector normalised such that its $L^2_\omega$-norm $||u_h^{(k)}||$ equals $1$. Since $u_h^{(k)}$ is orthogonal on $\ker\Delta_h$, it is also orthogonal on $\ker\Delta$ for all $h\neq 0$. For every $h>0$, we get

\begin{eqnarray}\label{eqn:delta_h-k_lbound}\delta_h^{(k)} = \langle\langle\Delta_hu_h^{(k)},\,u_h^{(k)}\rangle\rangle_\omega = \frac{h^2+1}{2}\,\langle\langle\Delta u_h^{(k)},\,u_h^{(k)}\rangle\rangle_\omega \geq \frac{h^2+1}{2}\,\delta^{(k)} \geq \frac{1}{2}\,\delta^{(k)},\end{eqnarray}

\noindent where $\delta^{(k)}>0$ is the smallest positive eigenvalue of $\Delta:C^\infty_k(X,\,\C)\longrightarrow C^\infty_k(X,\,\C)$ acting on $k$-forms. (So, $\delta^{(k)}$ is independent of $h$.)

Now, we know from Theorem 1.3 (and its corollary, Proposition 5.3) in [Pop17] that the Fr\"olicher spectral sequence of any compact Hermitian manifold $(X,\,\omega)$ degenerates at $E_1$ if and only if $\delta_h^{(k)}$ does not converge to zero at least as fast as $O(h^2)$ when $h\downarrow 0$ for every $k$. In our case, since the metric $\omega$ is K\"ahler, (\ref{eqn:delta_h-k_lbound}) shows that for every $k$, $\delta_h^{(k)}$ even remains uniformly bounded below by a positive constant when $h\downarrow 0$.  \hfill $\Box$

\vspace{3ex}

We can now infer the $d_h$-cohomology analogue of the standard $\partial\bar\partial$-lemma.

\begin{The}(the {\bf h-$\partial\bar\partial$-lemma})\label{The:h-dd-bar} Let $(X,\,\omega)$ be a compact {\bf K\"ahler} manifold with $\mbox{dim}_\C X=n$. For every $k\in\{0,1,\dots ,2n\}$, every $h\in\R\setminus\{0\}$ and every $k$-form $u\in\ker d_h\cap\ker d_{-h^{-1}}$, the following exactness conditions are equivalent:

$$u\in\mbox{Im}\,d_h \iff u\in\mbox{Im}\,d_{-h^{-1}} \iff u\in\mbox{Im}\,d \iff u\in\mbox{Im}\,(d_h\,d_{-h^{-1}}) = \mbox{Im}\,(\partial\bar\partial).$$

\end{The}

\noindent {\it Proof.} The equality $\mbox{Im}\,(d_h\,d_{-h^{-1}}) = \mbox{Im}\,(\partial\bar\partial)$ follows from (iii) of (\ref{eqn:obvious-prop}), while the property $u\in\mbox{Im}\,(d_h\,d_{-h^{-1}})$ obviously implies all the other exactness properties. 

Since $d_1=d$ and we allow any $h\neq 0$, it suffices to prove the implication ``$u\in\mbox{Im}\,d_h \implies  u\in\mbox{Im}\,(d_h\,d_{-h^{-1}})$`` for an arbitrary $h\neq 0$.  

Since $\Delta_h$ and $\Delta_{-h^{-1}}$ are self-adjoint elliptic operators with $d_h^2 = d_{-h^{-1}}^2=0$ and the manifold $X$ is compact, standard Hodge theory yields the following $L^2_\omega$-orthogonal decomposition (that does not require $\omega$ to be K\"ahler): \begin{eqnarray}\label{eqn:3-spqce-decomp} C^\infty_{k-1}(X,\,\C) = \ker\Delta_{-h^{-1}}\oplus\mbox{Im}\,d_{-h^{-1}}\oplus\mbox{Im}\,d_{-h^{-1}}^\star\end{eqnarray}

\noindent in which $\ker d_{-h^{-1}} = \ker\Delta_{-h^{-1}}\oplus\mbox{Im}\,d_{-h^{-1}}$.

Let $u\in C^\infty_k(X,\,\C)$ such that $u\in\ker d_h\cap\ker d_{-h^{-1}}$ and $u=d_h v$ for some $v\in C^\infty_{k-1}(X,\,\C)$. The $3$-space decomposition (\ref{eqn:3-spqce-decomp}) yields a unique decomposition 

$$v=v_0 + d_{-h^{-1}} u_1 + d_{-h^{-1}}^\star u_2,$$ 

\noindent where the $(k-1)$-form $v_0$ lies in $\ker\Delta_{-h^{-1}}$ and $u_1, u_2$ are of respective degrees $k-2$ and $k$. We get

$$u=d_h v = d_h v_0 + d_hd_{-h^{-1}} u_1 + d_hd_{-h^{-1}}^\star u_2 = -d_{-h^{-1}}d_h u_1 - d_{-h^{-1}}^\star d_h u_2.$$

\noindent Indeed, the last identity above follows from $v_0\in\ker\Delta_{-h^{-1}} = \ker\Delta_h = \ker d_h\cap\ker d_h^\star$ (where the K\"ahler assumption on $\omega$ was used to guarantee the proportionality of the Laplacians $\Delta_{-h^{-1}}$ and $\Delta_h$ -- see Theorem \ref{Thm:BKN} -- hence the equality of their kernels), from the anti-commutation of $d_h$ and $d_{-h^{-1}}$ (which holds trivially for any, not necessarily K\"ahler, metric $\omega$ -- see (iii) of (\ref{eqn:obvious-prop})) and from the anti-commutation of $d_h$ and $d_{-h^{-1}}^\star$ (which is a consequence of the K\"ahler assumption on $\omega$ via the $h$-commutation relations -- see Corollary \ref{Cor:Kaehler_anti-comm}).

Now, $u + d_{-h^{-1}}d_h u_1\in\ker d_{-h^{-1}}$ while $- d_{-h^{-1}}^\star d_hu_2\in\mbox{Im}\,d_{-h^{-1}}^\star$. However, $\ker d_{-h^{-1}}$ is orthogonal to $\mbox{Im}\,d_{-h^{-1}}^\star$, so the form $u + d_{-h^{-1}}d_h u_1 = - d_{-h^{-1}}^\star d_hu_2$, that lies in both subspaces, must vanish. In particular, $u = -d_{-h^{-1}}d_h u_1\in\mbox{Im}\,(d_h\,d_{-h^{-1}})$.  \hfill $\Box$

\vspace{3ex}

The above theorem leads naturally to the following

\begin{Def}\label{Def:h-dd-bar-manifolds} Let $h\in\R\setminus\{0\}$ be an arbitrary constant. A compact complex manifold $X$ with $\mbox{dim}_\C X=n$ is said to be an {\bf h-$\partial\bar\partial$-manifold} if for every $k\in\{0,1,\dots ,2n\}$ and every $k$-form $u\in\ker d_h\cap\ker d_{-h^{-1}}$, the following exactness conditions are equivalent:

$$u\in\mbox{Im}\,d_h \iff u\in\mbox{Im}\,d_{-h^{-1}} \iff u\in\mbox{Im}\,d \iff u\in\mbox{Im}\,(d_h\,d_{-h^{-1}}) = \mbox{Im}\,(\partial\bar\partial).$$

\end{Def}

Note that when $h=1$, $d_h=d$ and $d_{-h^{-1}}= d_{-1}$ coincides (up to a multiplicative constant) with $d^c$. The {\it h-$\partial\bar\partial$}-property introduced above does not require the form $u$ to be of pure type. In the cases $h\notin \{-1, 1\}$, it is meant to reinforce the standard $\partial\bar\partial$-property.

Like the standard $\partial\bar\partial$-property, the $h$-$\partial\bar\partial$-property is implied by the K\"ahler condition and implies the degeneration at the first page of the Fr\"oilicher spectral sequence (cf. Theorems \ref{The:h-dd-bar} above and \ref{The:h-dd-bar_E1} below). Actually, this last implication follows from the well-known implication with the $\partial\bar\partial$-property in place of the $h$-$\partial\bar\partial$-property, but we prefer to give a self-contained proof.

\begin{The}\label{The:h-dd-bar_E1} Let $h\in\R\setminus\{0\}$ be an arbitrary constant. The Fr\"olicher spectral sequence of any {\bf h-$\partial\bar\partial$-manifold} degenerates at $E_1$.

\end{The}

\noindent {\it Proof.} Let $X$ be an $h$-$\partial\bar\partial$-manifold with $\mbox{dim}_\C X=n$. For any bidegree $(p,\,q)$, pick any class $[\alpha]_{\bar\partial}\in E_1^{p,\,q}(X)$ and any representative $\alpha$ of $[\alpha]_{\bar\partial}$. We have $d_1([\alpha]_{\bar\partial}) = [\partial\alpha]_{\bar\partial}$.

 Moreover, since $\bar\partial\alpha=0$, we have $\partial\alpha = h\partial (h^{-1}\,\alpha) + \bar\partial (h^{-1}\,\alpha) = d_h(h^{-1}\,\alpha)\in\mbox{Im}\,d_h$. In particular, $\partial\alpha\in\ker d_h$ and $d_{-h^{-1}}(d_h(h^{-1}\,\alpha)) = ((h^2+1)/h^2)\,\partial\bar\partial\alpha = 0$, so $d_h(h^{-1}\,\alpha)\in\ker d_h\cap\ker d_{-h^{-1}}$. Thus, thanks to the $h$-$\partial\bar\partial$ assumption on $X$, the $d_h$-exactness of $\partial\alpha=d_h(h^{-1}\,\alpha)$ implies its $\partial\bar\partial$-exactness. In particular, $\partial\alpha\in\mbox{Im}\,\bar\partial$, hence $d_1([\alpha]_{\bar\partial}) = [\partial\alpha]_{\bar\partial} =0\in E_1^{p+1,\,q}(X)$.

This proves that all the differentials $d_1$ vanish identically, so $E_1^{p,\,q}(X) = E_2^{p,\,q}(X)$ for all $p,q$. 

Furthermore, since $\partial\alpha$ is $\partial\bar\partial$-exact, there exists a $(p,\,q-1)$-form $u$ such that $\partial\alpha = \bar\partial\partial u$, so $d_2([[\alpha]_{\bar\partial}]_{d_1}) = [[\partial(\partial u)]_{\bar\partial}]_{d_1} = 0\in E_2^{p+2,\,q-1}(X)$ and $d_r([\dots[[\alpha]_{\bar\partial}]_{d_1}\dots]_{d_{r-1}})= 0\in E_r^{p+r,\,q-r+1}(X)$ for all $r\geq 2$. 

Thus, all the differentials $d_r$ with $r\geq 1$ vanish identically. Hence, the Fr\"olicher spectral sequence of $X$ degenerates at $E_1$.  \hfill $\Box$

\subsection{The $h$-Bott-Chern and $h$-Aeppli cohomologies}\label{subsection:h-BC-A}

We start by defining the $h$-twisted analogues of the Bott-Chern and Aeppli cohomologies and by observing some basic properties of them. Unlike their standard counterparts, they are not defined in a given bidegree, but in a given total degree. 

\begin{Def}\label{Def:h-BC-A} Let $X$ be a compact complex $n$-dimensional manifold. For every $h\in\R\setminus\{0\}$ and every $k\in\{0,\dots , 2n\}$, we define the $k^{th}$ degree {\bf h-Bott-Chern} and {\bf h-Aeppli} cohomology groups by the formulae

$$H^k_{h-BC}(X,\,\C) = \frac{\ker d_h\cap\ker d_{-\frac{1}{h}}}{\mbox{Im}\,(d_hd_{-\frac{1}{h}})}  \hspace{5ex} \mbox{and} \hspace{5ex} H^k_{h-A}(X,\,\C) = \frac{\ker (d_hd_{-\frac{1}{h}})}{\mbox{Im}\,d_h + \mbox{Im}\,d_{-\frac{1}{h}}},$$

\noindent where all the vector spaces involved are subspaces of the space $C^\infty_k(X,\,\C)$ of smooth $k$-forms on $X$.

\end{Def}

We now observe some basic properties of these spaces that parallel their standard counterparts.

\begin{Lem}\label{Lem:h-dd-bar_basic-prop} Let $X$ be an $n$-dimensional compact complex manifold. 

\vspace{1ex}

(a)\, For every $h\in\R\setminus\{0\}$ and every $k\in\{0,\dots , 2n\}$, the canonical map $$T_h^{(k)}: H^k_{h-BC}(X,\,\C)\longrightarrow H^k_{h-A}(X,\,\C), \hspace{3ex} [\alpha]_{h-BC}\mapsto[\alpha]_{h-A},$$

\noindent is well defined. Moreover, if $X$ is an {\bf $h$-$\partial\bar\partial$-manifold} for some fixed $h\in\R\setminus\{0\}$, the map $T_h^{(k)}$ is an {\bf isomorphism} for every $k\in\{0,\dots , 2n\}$.

\vspace{1ex}

(b)\, For every $h\in\R\setminus\{0\}$ and every $k\in\{0,\dots , 2n\}$, the following identities hold: \begin{eqnarray}\label{eqn:h-classical_decomp}\nonumber H^k_{h-BC}(X,\,\C) & = & \bigoplus\limits_{p+q=k}H^{p,\,q}_{BC}(X,\,\C), \\
\nonumber H^k_{h-A}(X,\,\C) & = & \bigoplus\limits_{p+q=k}H^{p,\,q}_A(X,\,\C).\end{eqnarray}

\noindent Hence, the dimensions of the vector spaces $H^k_{h-BC}(X,\,\C)$ and $H^k_{h-A}(X,\,\C)$ are independent of $h$.

\vspace{1ex}

(c)\, For every $h\in\R\setminus\{0\}$ and every $k\in\{0,\dots , 2n\}$, the canonical maps

$$H^k_{h-BC}(X,\,\C)\longrightarrow H^k_{d_h}(X,\,\C)\longrightarrow H^k_{h-A}(X,\,\C), \hspace{3ex} [\alpha]_{h-BC}\mapsto [\alpha]_{d_h}\mapsto [\alpha]_{h-A},$$

\noindent are well defined. Moreover, if $X$ is an {\bf $h$-$\partial\bar\partial$-manifold} for some fixed $h\in\R\setminus\{0\}$, they are {\bf isomorphisms}, in particular their dimensions equal the $k^{th}$ Betti number $b_k$ of $X$, for every $k\in\{0,\dots , 2n\}$.

\end{Lem}

\noindent {\it Proof.} (a)\, Let $[\alpha]_{h-BC}\in H^k_{h-BC}(X,\,\C)$ be an arbitrary class and let $\alpha$ be an arbitrary representative of it. Then $d_h\alpha=0$ and $d_{-\frac{1}{h}}\alpha =0$, hence $d_hd_{-\frac{1}{h}}\alpha=0$, so $\alpha$ defines a class in $H^k_{h-A}(X,\,\C)$. To show that $[\alpha]_{h-A}$ is independent of the choice of representative $\alpha$ of the original class $[\alpha]_{h-BC}$, we have to show that $[\alpha]_{h-A}=0$ whenever $[\alpha]_{h-BC}=0$. However, this is obvious since $\mbox{Im}\,(d_hd_{-\frac{1}{h}})\subset\mbox{Im}\,d_h + \mbox{Im}\,d_{-\frac{1}{h}}$.

Suppose now that $X$ is an {\bf $h$-$\partial\bar\partial$-manifold} for some fixed $h\in\R\setminus\{0\}$. Fix any $k$. 

To show that $T_h^{(k)}$ is injective, suppose that $d_h\alpha=0$, $d_{-\frac{1}{h}}\alpha=0$ (i.e. $\alpha$ defines a class $[\alpha]_{h-BC}$) and $[\alpha]_{h-A}=0$ (i.e. $T_h^{(k)}([\alpha]_{h-BC})=0)$. In particular, $\alpha=d_h u + d_{-\frac{1}{h}}v$ for some forms $u,v$. Then $\alpha - d_h u = d_{-\frac{1}{h}}v\in\ker d_h\cap\mbox{Im}\,d_{-\frac{1}{h}}$, so $d_{-\frac{1}{h}}v\in\mbox{Im}\,(d_hd_{-\frac{1}{h}})$ thanks to the $h$-$\partial\bar\partial$-assumption. Meanwhile, $\alpha - d_{-\frac{1}{h}}v =d_h u\in\ker d_{-\frac{1}{h}}\cap\mbox{Im}\,d_h$, so $d_h u\in\mbox{Im}\,(d_hd_{-\frac{1}{h}})$ thanks to the $h$-$\partial\bar\partial$-assumption. Consequently, $\alpha=d_h u + d_{-\frac{1}{h}}v\in\mbox{Im}\,(d_hd_{-\frac{1}{h}})$, so $[\alpha]_{h-BC}=0$.

To show that $T_h^{(k)}$ is surjective, let $\alpha\in C^\infty_k(X,\,\C)$ such that $d_hd_{-\frac{1}{h}}\alpha=0$. We need to prove the existence of $(k-1)$-forms $u,v$ such that $d_h(\alpha + d_h u + d_{-\frac{1}{h}}v)=0$ and $d_{-\frac{1}{h}}(\alpha + d_h u + d_{-\frac{1}{h}}v)=0$. (Indeed, we will then have $[\alpha]_{h-A} = [\alpha + d_h u + d_{-\frac{1}{h}}v]_{h-A} = T_h^{(k)}([\alpha + d_h u + d_{-\frac{1}{h}}v]_{h-BC})$ with $[\alpha + d_h u + d_{-\frac{1}{h}}v]_{h-BC}$ well defined.) These identities are equivalent to $d_hd_{-\frac{1}{h}}v = -d_h\alpha$ and $d_{-\frac{1}{h}}d_h u = -d_{-\frac{1}{h}}\alpha$. Since $d_h\alpha\in\mbox{Im}\,d_h$ and $d_{-\frac{1}{h}}\alpha\in\mbox{Im}\,d_{-\frac{1}{h}}$ while both forms are simultaneously $d_h$-closed and $d_{-\frac{1}{h}}$-closed, they must be $(d_hd_{-\frac{1}{h}})$-exact thanks to the $h$-$\partial\bar\partial$-assumption. The surjectivity statement follows.

\vspace{1ex}

(b)\, For every $h\in\R\setminus\{0\}$ and every $k\in\{0,\dots , n\}$, the following equalities of subspaces of $C^\infty_k(X,\,\C)$ hold: \begin{eqnarray}\label{eqn:h-classic_eq_bis}\nonumber \ker d_h\cap\ker d_{-\frac{1}{h}} & = & \ker\partial\cap\ker\bar\partial  \\
 \mbox{Im}\,d_h + \mbox{Im}\,d_{-\frac{1}{h}} & = & \mbox{Im}\,\partial + \mbox{Im}\,\bar\partial.\end{eqnarray}

\noindent Indeed, for any $k$-form $\alpha$, the relation $\alpha\in\ker d_h\cap\ker d_{-\frac{1}{h}}$ is equivalent to having $h\partial\alpha + \bar\partial\alpha=0$ and $-\frac{1}{h}\,\partial\alpha + \bar\partial\alpha=0$, whose difference yields $(h+\frac{1}{h})\,\partial\alpha=0$, hence $\partial\alpha=0$ and $\bar\partial\alpha=0$. The reverse inclusion $\ker\partial\cap\ker\bar\partial\subset\ker d_h\cap\ker d_{-\frac{1}{h}}$ is obvious. Meanwhile, the relation $\alpha\in\mbox{Im}\,d_h + \mbox{Im}\,d_{-\frac{1}{h}}$ is equivalent to the existence of $(k-1)$-forms $u,v$ such that $\alpha = d_h u + d_{-\frac{1}{h}}v = \partial(h\,u - \frac{1}{h}\,v) + \bar\partial(u+v)$. 

Thus, for every $k$-form $\alpha = \sum\limits_{p+q=k}\alpha^{p,\,q}$ (written with its pure-type splitting apparent), the requirement $\alpha\in\ker d_h\cap\ker d_{-\frac{1}{h}}  =  \ker\partial\cap\ker\bar\partial$ is equivalent to the requirements $$\sum\limits_{p+q=k}\partial\alpha^{p,\,q}=0 \hspace{3ex} \mbox{and} \hspace{3ex}  \sum\limits_{p+q=k}\bar\partial\alpha^{p,\,q}=0,$$ 

\noindent which, in turn, are equivalent to requiring $\alpha^{p,\,q}\in\ker\partial\cap\ker\bar\partial = \ker d_h\cap\ker d_{-\frac{1}{h}}$ for all $p,q$. Similarly, thanks to (iii) of (\ref{eqn:obvious-prop}), requiring $\alpha\in\mbox{Im}\,(d_hd_{-\frac{1}{h}}) = \mbox{Im}(\partial\bar\partial)$ is equivalent to requiring $\alpha^{p,\,q}\in\mbox{Im}(\partial\bar\partial)$ for every $p,q$. We thus get the first decomposition of vector spaces stated under (b). The second decomposition is obtained in a similar way from the second identity in (\ref{eqn:h-classic_eq_bis}) and from (iii) of (\ref{eqn:obvious-prop}).

\vspace{1ex}

(c)\, The well-definedness of these maps follows at once from the inclusions $\ker d_h\cap\ker d_{-\frac{1}{h}}\subset\ker d_h\subset\ker\,(d_hd_{-\frac{1}{h}})$ and $\mbox{Im}\,(d_hd_{-\frac{1}{h}})\subset\mbox{Im}\,d_h\subset(\mbox{Im}\,d_h + \mbox{Im}\,d_{-\frac{1}{h}})$. The bijectivity of these maps when $X$ is supposed to be an $h$-$\partial\bar\partial$-manifold follows from straightforward applications of this hypothesis, from the proof of (a) and from the following lemma.  \hfill $\Box$

\begin{Lem}\label{Lem:h-dd-bar_closed-rep} If $X$ is an $h$-$\partial\bar\partial$-manifold, every $d_h$-cohomology class $[\alpha]_{d_h}$ (of any degree) contains a representative lying in $\ker d_h\cap\ker d_{-\frac{1}{h}}$.

\end{Lem}

\noindent {\it Proof.} Let $\alpha$ be a smooth $k$-form such that $d_h\alpha=0$. We wish to prove the existence of a smooth $(k-1)$-form $\beta$ such that $d_{-\frac{1}{h}}(\alpha + d_h\beta)=0$. This amounts to $d_{-\frac{1}{h}}d_h\beta = - d_{-\frac{1}{h}}\alpha$. However, $d_{-\frac{1}{h}}\alpha\in\ker d_h\cap\mbox{Im}\,d_{-\frac{1}{h}}$, so the $h$-$\partial\bar\partial$-hypothesis ensures that $d_{-\frac{1}{h}}\alpha\in\mbox{Im}\,(d_hd_{-\frac{1}{h}})$, proving the existence of $\beta$.  \hfill $\Box$

\vspace{3ex}

Recall that it was proved by Angella and Tomassini in [AT12] that on every compact complex manifold $X$, the inequality

\begin{equation}\label{eqn:AT_ineq} 2b_k\leq\sum\limits_{p+q=k}h^{p,\,q}_{BC} + \sum\limits_{p+q=k}h^{p,\,q}_A\end{equation}

\noindent holds for every $k$. Thanks to (b) of Lemma \ref{Lem:h-dd-bar_basic-prop}, this translates in our language to

\begin{equation}\label{eqn:AT_ineq_h-version} 2b_k\leq h^k_{h-BC} + h^k_{h-A}, \hspace{3ex} \mbox{for all} \hspace{1ex} k\in\{0,\dots , 2n\} \hspace{1ex} \mbox{and all} \hspace{1ex} h\in\R\setminus\{0\},\end{equation}

\noindent where $h^k_{h-BC}:=\mbox{dim} H^k_{h-BC}(X,\,\C)$ and $h^k_{h-A}:=\mbox{dim} H^k_{h-A}(X,\,\C)$. (Recall that we always have $b_k= \mbox{dim} H^k_{d_h}(X,\,\C)$ for all $k$ and $h\neq 0$, see e.g. Introduction.) Moreover, the second main result of [AT12] states that equality holds in (\ref{eqn:AT_ineq}) for every $k$ if and only if $X$ satisfies {\it version (a) of the $\partial\bar\partial$-lemma} (see Introduction).

\begin{Cor}\label{Cor:h-dd-bar_non-jumping} (a)\footnote{This is already obvious from the definitions, but we give a new argument to show the consistency of several results with one another.}\, Let $X$ be an $n$-dimensional compact complex manifold. If $X$ is an $h$-$\partial\bar\partial$-manifold for some $h\in\R\setminus\{0\}$, then $X$ satisfies version (a) of the $\partial\bar\partial$-lemma (see Introduction).

\vspace{1ex}

(b)\, Let $(X_t)_{t\in\Delta}$ be a holomorphic family of compact complex manifolds. If some fibre $X_0$ is an $h$-$\partial\bar\partial$-manifold for some $h\in\R\setminus\{0\}$, then the $h$-Bott-Chern numbers $h^k_{h-BC}(t):=\mbox{dim} H^k_{h-BC}(X_t,\,\C) $ and the $h$-Aeppli numbers $h^k_{h-A}(t):=\mbox{dim} H^k_{h-A}(X_t,\,\C)$ remain constant in a neighbourhood of $X_0$:

$$h^k_{h-BC}(t) = h^k_{h-BC}(0) \hspace{3ex} \mbox{and} \hspace{3ex} h^k_{h-A}(t)=h^k_{h-A}(0)$$

\noindent for all $k\in\{0,\dots , 2n\}$ and all $t\in\Delta$ close enough to $0$.

\end{Cor}

\noindent {\it Proof.} (a)\, If $X$ is an $h$-$\partial\bar\partial$-manifold for some $h\in\R\setminus\{0\}$, it follows from (c) of Lemma \ref{Lem:h-dd-bar_basic-prop} that $2b_k = h^k_{h-BC} + h^k_{h-A}$ for every $k$. This is equivalent to $X$ satisfying {\it version (a) of the $\partial\bar\partial$-lemma} by the above discussion and the second main result of [AT12].

(b)\, It follows from (b) of Lemma \ref{Lem:h-dd-bar_basic-prop} that $h^k_{h-BC}(t) = \sum\limits_{p+q=k}h^{p,\,q}_{BC}(t)$ and $h^k_{h-A}(t) = \sum\limits_{p+q=k}h^{p,\,q}_A(t)$ for all $k$ and all $t$. Since the Bott-Chern and the Aeppli numbers are known to satisfy the semicontinuity property $h^{p,\,q}_{BC}(0)\geq h^{p,\,q}_{BC}(t)$ and $h^{p,\,q}_A(0)\geq h^{p,\,q}_A(t)$ for all $t$ close enough to $0$ and all $p,q$, we infer the analogous property for the $h$-Bott-Chern and the $h$-Aeppli numbers:

\begin{equation}\label{eqn:semicontinuity_h-BC-A}h^k_{h-BC}(0)\geq h^k_{h-BC}(t) \hspace{3ex} \mbox{and} \hspace{3ex} h^k_{h-A}(0)\geq h^k_{h-A}(t)\end{equation}

\noindent for all $t$ close enough to $0$ and all $k$. 

 On the other hand, since $X_0$ is an $h$-$\partial\bar\partial$-manifold for some $h\in\R\setminus\{0\}$, we have $2b_k = h^k_{h-BC}(0) + h^k_{h-A}(0)$ for every $k$ (cf. (c) of Lemma \ref{Lem:h-dd-bar_basic-prop}). Hence, $2b_k \geq h^k_{h-BC}(t) + h^k_{h-A}(t)$ for every $k$ and every $t$ close enough to $0$ thanks also to (\ref{eqn:semicontinuity_h-BC-A}). However, the reverse inequality  $2b_k \leq h^k_{h-BC}(t) + h^k_{h-A}(t)$ also holds for all $t$ and $k$ thanks to [AT12] (cf. (\ref{eqn:AT_ineq_h-version})). The result follows.  \hfill $\Box$

\subsection{Deformation openness of the $h$-$\partial\bar\partial$-property}\label{subsection:h-openness}

We now prove the following analogue for $h$-$\partial\bar\partial$-manifolds of Wu's openness result for $\partial\bar\partial$-manifolds (cf. [Wu06]).

\begin{The}\label{The:h-dd-bar-openness} Let $\pi:{\cal X}\longrightarrow\Delta$ be a proper holomorphic submersion from a complex manifold ${\cal X}$ to a ball $\Delta\subset\C^N$ containing the origin. For every $t\in\Delta$, let $X_t:=\pi^{-1}(t)$ be the fibre above $t$. Fix an arbitrary constant $h\in\R\setminus\{0\}$.

If $X_0$ is an $h$-$\partial\bar\partial$-manifold, then $X_t$ is an $h$-$\partial\bar\partial$-manifold for all $t\in\Delta$ sufficiently close to $0$.

\end{The}

The proof will follow the pattern of the one given by Wu in [Wu06] for the deformation openness of the standard $\partial\bar\partial$-property. For the sake of consistency, we will follow the presentation in $\S.4.3$ of [Pop14] where Wu's arguments and some of those in [DGMS75] were re-expalined, while pointing out the changes needed in our current $h$-$\partial\bar\partial$-context.  

We start with some ad hoc terminology that parallels Definition 4.7. in [Pop14].

\begin{Def}\label{Def:prop} Let $h\in\R\setminus\{0\}$ be an arbitrary constant. For any given $k=0, 1, \dots , 2n$, a given $n$-dimensional compact complex manifold $X$ is said to satisfy property: \\

 $(A_k)$ if the canonical map $H^k_{h-BC}(X, \, \C) \longrightarrow  H^k_{h-A}(X, \, \C)$ is {\bf injective}. \\

\noindent This property is equivalent to the property

\vspace{1ex}

$(A_k')$\,\, $\ker d_h\cap\ker d_{-\frac{1}{h}}\cap(\mbox{Im}\,d_h + \mbox{Im}\,d_{-\frac{1}{h}}) =\mbox{Im}\,(d_h\,d_{-\frac{1}{h}}) \hspace{2ex} \mbox{as subspaces of}\,\,C^\infty_k(X,\,\C).$ \\

 $(B_k)$ if the canonical map $H^k_{h-BC}(X, \, \C) \longrightarrow  H^k_{h-A}(X, \, \C)$ is {\bf surjective}. \\

\noindent This property is equivalent to the property

\vspace{1ex}

$(B_k')$\, $\mbox{Im}\,d_h + \mbox{Im}\,d_{-\frac{1}{h}} + (\ker d_h \cap \ker d_{-\frac{1}{h}})=\ker(d_h\,d_{-\frac{1}{h}}) \hspace{2ex} \mbox{as subspaces of}\,\, C^\infty_k(X,\,\C).$ \\

 $(C_k)$ if the canonical maps $H^k_{h-BC}(X, \, \C) \longrightarrow  H^k_{d_{-\frac{1}{h}}}(X, \, \C)$ and $H^k_{h-BC}(X, \, \C) \longrightarrow  H^k_{d_h}(X, \, \C)$ are {\bf injective}. \\

\noindent This property is equivalent to the simultaneous occurence of 

\vspace{1ex}

$(C_k')$$(i)$\, $\mbox{Im}\,d_{-\frac{1}{h}}\cap\ker d_h=\mbox{Im}\,(d_h\,d_{-\frac{1}{h}})$ \hspace{2ex} \mbox{and} \hspace{2ex} $(C_k')$$(ii)$\, $\mbox{Im}\,d_h\cap\ker d_{-\frac{1}{h}} = \mbox{Im}\,(d_h\,d_{-\frac{1}{h}})$

\vspace{2ex}

\noindent as subspaces of $C^\infty_k(X,\,\C)$. \\

$(D_k')$ if $(i)$\,  $\mbox{Im}\,d_h + \ker d_{-\frac{1}{h}} = \ker(d_h\,d_{-\frac{1}{h}})$ \hspace{1ex} and \hspace{1ex} $(ii)$\, $\mbox{Im}\,d_{-\frac{1}{h}} + \ker d_h = \ker(d_h\,d_{-\frac{1}{h}})$

\vspace{2ex}

\noindent as subspaces of $C^\infty_k(X,\,\C)$. \\

$(L_k)$ if for every $k$-form $u\in\ker d_h\cap\ker d_{-h^{-1}}$, the following exactness conditions are equivalent: $$u\in\mbox{Im}\,d_h \iff u\in\mbox{Im}\,d_{-h^{-1}} \iff u\in\mbox{Im}\,(d_h\,d_{-h^{-1}}) = \mbox{Im}\,(\partial\bar\partial).$$

\end{Def}

Property $(L_k)$ is a restatement of the pair of properties $(C'_k)(i)$ and $(C'_k)(ii)$. As already pointed out, the following equivalences are immediate: 

$$(A_k)\Longleftrightarrow (A_k'), \hspace{2ex} (B_k)\Longleftrightarrow (B_k'), \hspace{2ex} (C_k)\Longleftrightarrow (C_k').$$ 

\noindent Meanwhile, the inclusions $\supset$ in $(A_k')$, $\subset$ in $(B_k')$, $\supset$ in $(C_k')(i), (ii)$ and $\subset$ in $(D_k')(i), (ii)$ always hold trivially. The following statement is the $h$-$\partial\bar\partial$ analogue of a fact implicitly proved in [DGMS75] in the standard $\partial\bar\partial$ context and will provide a key ingredient for the proof of Theorem \ref{The:h-dd-bar-openness}.

\begin{Prop}(the $h$-$\partial\bar\partial$-analogue of Lemma 5.15 in [DGMS75])\label{Prop:h_DGMS-ddbar-char} Let $h\in\R\setminus\{0\}$ be an arbitrary constant. Let $X$ be a compact $n$-dimensional complex manifold. For every $k=1, \dots , 2n$, the following equivalences hold\!\!:

$$(L_k)\Longleftrightarrow (A_k)\Longleftrightarrow (C_k) \Longleftrightarrow (D_{k-1}') \Longleftrightarrow (B_{k-1}).$$

\end{Prop}

\noindent {\it Proof.} Fix an arbitrary $k\in\{1, \dots , 2n\}$. Given the above explanations, it suffices to prove the equivalences

$$(A_k')\Longleftrightarrow (C_k') \Longleftrightarrow (D_{k-1}') \Longleftrightarrow (B_{k-1}').$$

\noindent {\it Proof of $(A_k')\Longrightarrow (C_k').$} Let $u\in C^{\infty}_k(X, \, \C)$ such that $d_h u=0$ and $u=d_{-\frac{1}{h}} v$ for some $(k-1)$-form $v$. Then $u\in\ker d_h\cap\ker d_{-\frac{1}{h}}\cap(\mbox{Im}\,d_h + \mbox{Im}\,d_{-\frac{1}{h}})$. So $(A_k')$ forces $u\in\mbox{Im}\,(d_h\,d_{-\frac{1}{h}})$. This proves $(i)$ of $(C_k')$. The proof of $(ii)$ of $(C_k')$ is similar with $d_h$ and $d_{-\frac{1}{h}}$ reversed. \\

\noindent {\it Proof of $(C_k')\Longrightarrow (A_k').$} Let $u\in C^{\infty}_k(X, \, \C)$ such that $d_h u=0$, $d_{-\frac{1}{h}} u=0$ and $u=d_h v + d_{-\frac{1}{h}} w$ for some $(k-1)$-forms $v$ and $w$. Then  \\

$\cdot$\, $\mbox{Im}\,d_h\ni d_h v = u - d_{-\frac{1}{h}} w\in\ker d_{-\frac{1}{h}}$, so $d_h v\in \mbox{Im}\,d_h\cap\ker d_{-\frac{1}{h}} = \mbox{Im}\,(d_hd_{-\frac{1}{h}})$, the last identity of subspaces being given by the hypothesis $(C_k')(i)$.  \\

$\cdot$\, $\mbox{Im}\,d_{-\frac{1}{h}}\ni d_{-\frac{1}{h}} w = u - d_h v\in\ker d_h$, so $d_{-\frac{1}{h}} w\in \mbox{Im}\,d_{-\frac{1}{h}}\cap\ker d_h = \mbox{Im}\,(d_hd_{-\frac{1}{h}})$, the last identity of subspaces being given by the hypothesis $(C_k')(ii)$.  \\

 We now get $u=d_h v + d_{-\frac{1}{h}} w\in\mbox{Im}\,(d_hd_{-\frac{1}{h}})$. This proves $(A_k')$. \\

\noindent {\it Proof of $(C_k')\Longrightarrow (D_{k-1}').$} Let $u\in C^{\infty}_{k-1}(X, \, \C)$ such that $d_hd_{-\frac{1}{h}} u=0$. Then: \\

$\cdot$\, $d_{-\frac{1}{h}} u$ is a $k$-form and $d_{-\frac{1}{h}} u\in\ker d_h\cap\mbox{Im}\,d_{-\frac{1}{h}}=\mbox{Im}\,(d_hd_{-\frac{1}{h}})$, the last identity of subspaces being given by the hypothesis $(C_k')(i)$. So $d_{-\frac{1}{h}} u = d_{-\frac{1}{h}}d_h\zeta$ for some $(k-2)$-form $\zeta$. This amounts to $u-d_h\zeta\in\ker d_{-\frac{1}{h}}$. 

 We get $u=d_h\zeta + (u-d_h\zeta)\in\mbox{Im}\,d_h + \ker d_{-\frac{1}{h}}$. This proves $(D_{k-1}')(i).$ \\

$\cdot$\, $d_h u$ is a $k$-form and $d_h u\in\ker d_{-\frac{1}{h}}\cap\mbox{Im}\,d_h=\mbox{Im}\,(d_hd_{-\frac{1}{h}})$, the last identity of subspaces being given by the hypothesis $(C_k')(ii)$. Hence $d_h u = d_hd_{-\frac{1}{h}} w$ for some $(k-2)$-form $w$. This amounts to $u - d_{-\frac{1}{h}} w\in\ker d_h$.

 We get $u=d_{-\frac{1}{h}} w + (u - d_{-\frac{1}{h}} w)\in \mbox{Im}\,d_{-\frac{1}{h}} + \ker d_h$. This proves $(D_{k-1}')(ii).$ \\

\noindent {\it Proof of $(D_{k-1}')\Longrightarrow (C_k').$} Let $u\in C^{\infty}_k(X, \, \C)$ such that $d_h u=0$ and $u=d_{-\frac{1}{h}} v$ for some $(k-1)$-form $v$. Then $v\in\ker(d_hd_{-\frac{1}{h}})= \mbox{Im}\,d_h + \ker d_{-\frac{1}{h}}$, where the last identity of subspaces is $(D_{k-1}')(i)$. Thus, we can find a $(k-2)$-form $w$ and a $(k-1)$-form $\zeta$ such that 

$$v=d_h w + \zeta \hspace{2ex} \mbox{and} \hspace{2ex} d_{-\frac{1}{h}}\zeta=0.$$

\noindent Applying $d_{-\frac{1}{h}}$, we get: $u=d_{-\frac{1}{h}} v = d_{-\frac{1}{h}}d_h w\in\mbox{Im}\,(d_hd_{-\frac{1}{h}})$. This proves $(C_k')(i)$.

 Reversing the roles of $d_h$ and $d_{-\frac{1}{h}}$, we get $(C_k')(ii)$ in a similar way from $(D_{k-1}')(ii)$. \\

\noindent {\it Proof of $(D_{k-1}')\Longrightarrow (B_{k-1}').$} Let $u\in C^{\infty}_{k-1}(X, \, \C)$ such that $d_hd_{-\frac{1}{h}} u=0$. Thanks to $(D_{k-1}')(ii)$, we can find a $(k-2)$-form $v$ and a $(k-1)$-form $w$ such that 

$$u = d_{-\frac{1}{h}} v + w  \hspace{2ex} \mbox{and} \hspace{2ex} w\in\ker d_h.$$

\noindent Thus $d_hd_{-\frac{1}{h}} w=0$, so by $(D_{k-1}')(i)$ we can write

$$w = d_h\zeta + \rho  \hspace{2ex} \mbox{with} \hspace{2ex} \rho\in\ker d_{-\frac{1}{h}}$$

 \noindent for some $(k-2)$-form $\zeta$ and some $(k-1)$-form $\rho$. We get $\rho = w - d_h\zeta\in\ker d_h$ (because $w\in\ker d_h$). Given the choice of $\rho$, this implies that $\rho\in\ker d_h\cap\ker d_{-\frac{1}{h}}$.  

 Putting the bits together, we have

$$u = d_{-\frac{1}{h}} v + d_h\zeta + \rho\in\mbox{Im}\,d_{-\frac{1}{h}} + \mbox{Im}\,d_h + (\ker d_h\cap\ker d_{-\frac{1}{h}}).$$ This proves $(B_{k-1}')$. \\

\noindent {\it Proof of $(B_{k-1}')\Longrightarrow (D_{k-1}').$} This implication is trivial because $\mbox{Im}\,d_{-\frac{1}{h}} + (\ker d_h\cap\ker d_{-\frac{1}{h}}) \subset\ker d_{-\frac{1}{h}}$ and $\mbox{Im}\,d_h + (\ker d_h\cap\ker d_{-\frac{1}{h}}) \subset\ker d_h$. \\

 The proof of Proposition \ref{Prop:h_DGMS-ddbar-char} is complete.  \hfill $\Box$

\vspace{3ex}

Note that the simultaneous occurence of properties $(L_k)$ for all $k\in\{0,\dots , 2n\}$ is an a priori weaker condition than the $h$-$\partial\bar\partial$-property since it does not include the equivalence $u\in\mbox{Im}\,d \iff u\in\mbox{Im}\,(d_h\,d_{-h^{-1}})$. However, we can easily see as a consequence of Proposition \ref{Prop:h_DGMS-ddbar-char} that these two conditions are actually equivalent.

\begin{Cor}\label{Cor:equivalence_h-dd-bar_properties} Let $h\in\R\setminus\{0\}$ be an arbitrary constant. Let $X$ be a compact complex manifold with $\mbox{dim}_\C X=n$. Fix an arbitrary $k\in\{0,\dots ,2n\}$ and suppose that for every $k$-form $u\in\ker d_h\cap\ker d_{-h^{-1}}$, the following exactness conditions are equivalent: $$u\in\mbox{Im}\,d_h \iff u\in\mbox{Im}\,d_{-h^{-1}} \iff u\in\mbox{Im}\,(d_h\,d_{-h^{-1}}) = \mbox{Im}\,(\partial\bar\partial).$$

Then, for every $k$-form $u\in\ker d_h\cap\ker d_{-h^{-1}}$, the equivalence ``$u\in\mbox{Im}\,d \iff u\in\mbox{Im}\,(d_h\,d_{-h^{-1}}) = \mbox{Im}\,(\partial\bar\partial)$'' also holds.

In particular, if the assumption is made for all $k\in\{0,\dots ,2n\}$, then $X$ is an $h$-$\partial\bar\partial$-manifold.

\end{Cor}

\noindent {\it Proof.} Since the implication $u\in\mbox{Im}\,(d_h\,d_{-h^{-1}}) \Longrightarrow u\in\mbox{Im}\,d$ is trivial, we only have to prove the reverse implication for every $k$-form $u\in\ker d_h\cap\ker d_{-h^{-1}}$. Let $u\in\mbox{Im}\,d$ be such a $k$-form. Then $u\in\mbox{Im}\,d_h + \mbox{Im}\,d_{-h^{-1}}$ thanks to identity (iv) in (\ref{eqn:obvious-prop}), so $u$ defines a class $[u]_{h-BC}\in H^k_{h-BC}(X,\,\C)$ that maps to the zero class in $H^k_{h-A}(X,\,\C)$ under the canonical map $H^k_{h-BC}(X, \, \C) \longrightarrow  H^k_{h-A}(X, \, \C)$. Then by $(A_k)$, which holds because it is equivalent to our assumption $(L_k)$ thanks to Proposition \ref{Prop:h_DGMS-ddbar-char}, this map is injective. Hence, $[u]_{h-BC}=0\in H^k_{h-BC}(X,\,\C)$, so $u\in\mbox{Im}\,(d_h\,d_{-h^{-1}})$.  \hfill $\Box$

\vspace{3ex}

We are now well equipped to prove the deformation openness of the $h$-$\partial\bar\partial$-property of compact complex manifolds.

\vspace{3ex}

\noindent {\it Proof of Theorem \ref{The:h-dd-bar-openness}.} The arguments are analogues in the $h$-$\partial\bar\partial$ context of those given by Wu in the classical $\partial\bar\partial$ context. As in [Wu06], the main idea is to exploit, for every fixed $k$, the equivalence $$(L_k) \Longleftrightarrow (A_k) \Longleftrightarrow (B_{k-1}),$$ \noindent namely the discrepancy of one degree between the characterisation of the $h$-$\partial\bar\partial$-property for $k$-forms in terms of the injectivity of the $h$-BC$\to h$-A-map and in terms of its surjectivity. This prompts an argument by induction on $k$, since the $h$-$\partial\bar\partial$-property holds trivially in degree $k=0$ (i.e. for functions). 

To show that $X_t$ is an $h$-$\partial\bar\partial$-manifold for all $t\in\Delta$ sufficiently close to $0$, suppose that the $h$-$\partial\bar\partial$-property holds in degree $k$ on $X_t$ for the operators $d_h(t)$ and $d_{-\frac{1}{h}}(t)$ induced by the complex structure of $X_t$ for all $t$ close to $0$. We will prove that the same is true in degree $k+1$. 

 The $h$-$\partial\bar\partial$ assumption on $X_0$ in degree $k$ implies the following identities thanks to (c) of Lemma \ref{Lem:h-dd-bar_basic-prop} and respectively (b) of Corollary \ref{Cor:h-dd-bar_non-jumping}: \begin{eqnarray}\nonumber\mbox{dim}_{\C}H^k_{h-BC}(X_0, \, \C) & = & \mbox{dim}_{\C}H^k_{h-A}(X_0, \, \C),\\
\nonumber \mbox{dim}_{\C}H^k_{h-BC}(X_0, \, \C) & = & \mbox{dim}_{\C}H^k_{h-BC}(X_t, \, \C) \hspace{6ex} \mbox{for all} \hspace{1ex} t \hspace{1ex} \mbox{close to} \hspace{1ex} 0,\\
\nonumber \mbox{dim}_{\C}H^k_{h-A}(X_0, \, \C) & = & \mbox{dim}_{\C}H^k_{h-A}(X_t, \, \C) \hspace{6ex} \mbox{for all} \hspace{1ex} t \hspace{1ex} \mbox{close to} \hspace{1ex} 0.\end{eqnarray}

\noindent Hence, $\mbox{dim}_{\C}H^k_{h-BC}(X_t, \, \C) = \mbox{dim}_{\C}H^k_{h-A}(X_t, \, \C)$ for all $t\in\Delta$ close to $0$.

Meanwhile, by Proposition \ref{Prop:h_DGMS-ddbar-char}, the induction hypothesis $(L_k)$ on $X_t$ is equivalent to the canonical linear map $H^k_{h-BC}(X_t,\,\C)\to H^k_{h-A}(X_t,\,\C)$ being {\it injective} (property $(A_k)$). Since these are {\it finite-dimensional} vector spaces of {\it equal dimensions}, the linear map $H^k_{h-BC}(X_t,\,\C)\to H^k_{h-A}(X_t,\,\C)$ must also be {\it surjective}. Thus, property $(B_k)$ holds on $X_t$ for all $t\in\Delta$ close to $0$. However, thanks to Proposition \ref{Prop:h_DGMS-ddbar-char}, this is equivalent to property $(L_{k+1})$, i.e. to the $h$-$\partial\bar\partial$-property in degree $k+1$, holding on $X_t$ for all $t\in\Delta$ close to $0$.   \hfill $\Box$

\vspace{6ex}

\noindent {\bf References.} \\

\noindent [AT12]\, D. Angella, A. Tomassini --- {\it On the $\partial\bar\partial$-Lemma and Bott-Chern Cohomology} --- Invent. Math.

\vspace{1ex}

\noindent [COUV16]\, M. Ceballos, A. Otal, L. Ugarte, R. Villacampa --- {\it Invariant Complex Structures on $6$-Nilmanifolds: Classification, Fr\"olicher Spectral Sequence and Special Hermitian Metrics} --- J. Geom. Anal. {\bf 26} (2016), no. 1, 252–286.

\vspace{1ex}

\noindent [Dem84]\, J.-P. Demailly --- {\it Sur l'identit\'e de Bochner-Kodaira-Nakano en g\'eom\'etrie hermitienne} --- S\'eminaire d'analyse P. Lelong, P. Dolbeault, H. Skoda (editors) 1983/1984, Lecture Notes in Math., no. {\bf 1198}, Springer Verlag (1986), 88-97.

\vspace{1ex}

\noindent [Dem92]\, J.-P. Demailly --- {\it Regularization of Closed Positive Currents and Intersection Theory} --- J. Alg. Geom., {\bf 1} (1992), 361-409.

\vspace{1ex}

\noindent [DGMS75]\, P. Deligne, Ph. Griffiths, J. Morgan, D. Sullivan --- {\it Real Homotopy Theory of K\"ahler Manifolds} --- Invent. Math. {\bf 29} (1975), 245-274.

\vspace{1ex}

\noindent [Gau77]\, P. Gauduchon --- {\it Le th\'eor\`eme de l'excentricit\'e nulle} --- C.R. Acad. Sc. Paris, S\'erie A, t. {\bf 285} (1977), 387-390.

\vspace{1ex}

\noindent [KS60]\, K. Kodaira, D.C. Spencer --- {\it On Deformations of Complex Analytic Structures, III. Stability Theorems for Complex Structures} --- Ann. Math. {\bf 71}, No. 1 (1960), 43-76.

\vspace{1ex}

\noindent [Lam99]\, A. Lamari --- {Courants k\"ahl\'eriens et surfaces compactes} --- Ann. Inst. Fourier, Grenoble, {\bf 49}, 1 (1999), 263-285.

\vspace{1ex}

\noindent [Pop13]\, D. Popovici --- {\it Deformation Limits of Projective Manifolds\!\!: Hodge Numbers and Strongly Gauduchon Metrics} --- Invent. Math. {\bf 194} (2013), 515-534.

\vspace{1ex}

\noindent [Pop14]\, D. Popovici --- {\it Deformation Openness and Closedness of Various Classes of Compact Complex Manifolds; Examples} --- Ann. Sc. Norm. Super. Pisa Cl. Sci. (5), Vol. XIII (2014), 255-305. 

\vspace{1ex}

\noindent [Pop15]\, D. Popovici --- {\it Aeppli Cohomology Classes Associated with Gauduchon Metrics on Compact Complex Manifolds} --- Bull. Soc. Math. France {\bf 143} (3), (2015), p. 1-37.

\vspace{1ex}

\noindent [Pop16]\, D. Popovici --- {\it Degeneration at $E_2$ of Certain Spectral Sequences} --- International Journal of Mathematics {\bf 27}, no. 14 (2016), 1650111, 31 pp.

\vspace{1ex}

\noindent [Pop17]\, D. Popovici --- {\it Adiabatic Limit and the Fr\"olicher Spectral Sequence} --- arXiv e-print CV 1709.04332v1

\vspace{1ex}

\noindent [PU18a]\, D. Popovici, L. Ugarte --- {\it Compact Complex Manifolds with Small Gauduchon Cone} --- Proceedings of the London Mathematical Society (3) (2018) doi:10.1112/plms.12110.

\vspace{1ex}

\noindent [PU18b]\, D. Popovici, L. Ugarte --- {\it Symmetry and Duality for a $5$-Dimensional Nilmanifold} --- in preparation.

\vspace{1ex}

\noindent [Wu06]\, C.-C. Wu --- {\it On the Geometry of Superstrings with Torsion} --- thesis, Department of Mathematics, Harvard University, April 2006.

\vspace{6ex}

\noindent Ibn Tofail University, Faculty of Sciences, Departement of Mathematics, PO 242 Kenitra, Morocco

\noindent Email: houda.bellitir@uit.ac.ma

\vspace{2ex}

\noindent and

\vspace{2ex}

\noindent Universit\'e Paul Sabatier, Institut de Math\'ematiques de Toulouse

\noindent 118, route de Narbonne, 31062, Toulouse Cedex 9, France

\noindent Email: popovici@math.univ-toulouse.fr

\end{document}